\documentclass[11pt,reqno]{amsart}
\usepackage{amssymb}
\usepackage{amsfonts}
\usepackage{amsmath}
\usepackage{graphicx}
\usepackage{hyperref}
\usepackage{float}
\usepackage{epstopdf}
\usepackage{color}
\usepackage{bm}
\usepackage{comment}
\usepackage{soul}

\allowdisplaybreaks
\newtheorem{theorem}{Theorem}[section]
\newtheorem{proposition}[theorem]{Proposition}
\newtheorem{lemma}[theorem]{Lemma}
\newtheorem{corollary}[theorem]{Corollary}

\newtheorem{definition}[theorem]{Definition}

\newtheorem{thm}{Theorem}

\def\SG{\mathcal{SG}}

\def\mcE{\mathcal{E}}

\numberwithin{equation}{section}

\begin{document}
\title[Sobolev spaces on p.c.f. self-similar sets: critical orders and atomic decompositions]{Sobolev spaces on p.c.f. self-similar sets: critical orders and atomic decompositions}

\author{Shiping Cao}
\address{Department of Mathematics, Cornell University, Ithaca 14853, USA}
\email{sc2873@cornell.edu}
\thanks{}

\author{Hua Qiu}
\address{Department of Mathematics, Nanjing University, Nanjing 210093, China}
\email{huaqiu@nju.edu.cn}
\thanks{The research of Qiu was supported by the NSFC grant 11471157}


\date{}


\begin{abstract}
We consider the Sobolev type spaces $H^\sigma(K)$ with $\sigma\geq 0$, where $K$ is a post-critically finite self-similar set with the natural boundary. Firstly, we compare different classes of Sobolev spaces $H^\sigma_N(K),H^\sigma_D(K)$ and $H^\sigma(K)$, {and observe} a sequence of critical orders of $\sigma$ in our comparison theorem. Secondly, We present a general atomic decomposition theorem of Sobolev spaces $H^\sigma(K)$, where the same critical orders play an important role. At the same time, we provide purely analytic approaches for various Besov type characterizations of Sobolev spaces $H^\sigma(K)$.
\end{abstract}
\maketitle

\tableofcontents

\section{Introduction}\label{intro}
\begin{figure}[h]
	\includegraphics[width=5cm]{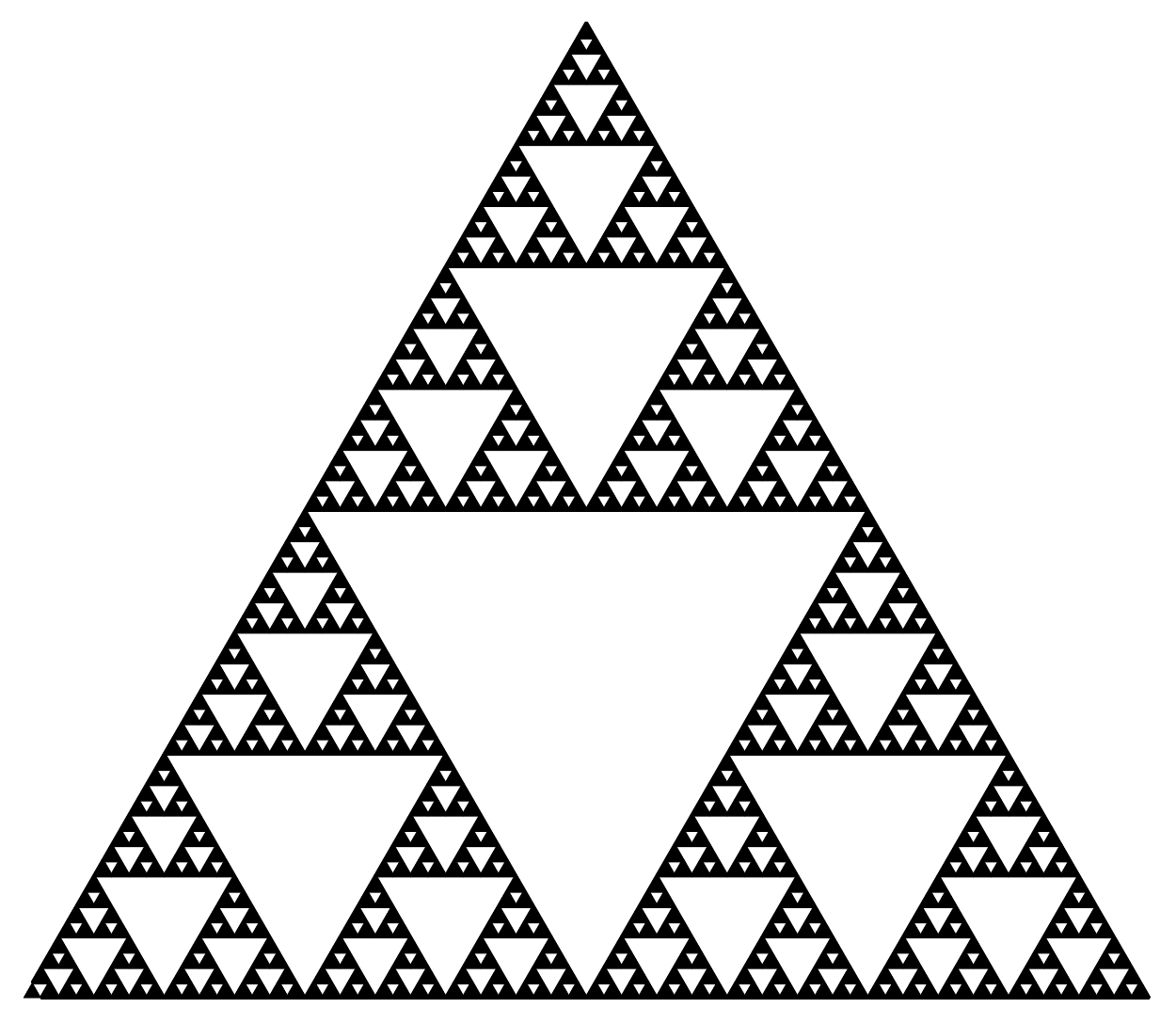}
	\begin{picture}(0,0)
	\put(-5,0){$p_3$}
	\put(-152,0){$p_2$}
	\put(-79,125){$p_1$}
	\end{picture}	
	\caption{the Sierpinski gasket.}\label{SG}
\end{figure}
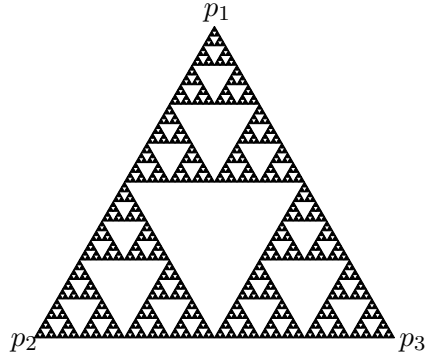

Analysis on fractals, based on the construction of Laplacians and Dirichlet forms, has been studied for years. In this paper, we study function spaces on post-critically finite (p.c.f.) self-similar sets, based on Kigami's construction of Dirichlet forms (\cite{k2,k3}). See backgrounds in books \cite{book1,book2}, and see \cite{Ba} for the probabilisitic approach. One of the most well-known example of p.c.f. self-similar sets is the Sierpinski gasket ($\SG$), see Figure \ref{SG}.

In 2003, on a p.c.f. fractal $K$, Strichartz\cite{Functionspace} considered general Sobolev type spaces $L^p_\sigma(K)$ and Besov type spaces $\Lambda^{p,q}_\sigma(K)$ with $1\leq p,q\leq \infty$ and $\sigma\geq 0$, and gave an all-round study on various embedding theorems and interpolation results. A systematical introduction of Sobolev spaces can also be found in \cite{Gr} on more general metric measure spaces by Grigor'yan. We will focus on the important case $p=q=2$ in this paper, and use the notation $H^\sigma(K)$ for these $L^2$ Sobolev spaces. In the definition, the p.c.f. self-similar $K$ is viewed as a domain with the natural boundary, which consists of finitely many points. For example, the Sierpinski gasket $\mathcal{SG}$ has the boundary $V_0=\{p_1,p_2,p_3\}$, the three vertices of a triangle, as shown in Figure \ref{SG}.

On the other hand, one may define the Sobolev spaces $H^\sigma_N(K)$ with the Bessel potential {$(1-\Delta_N)^{-\sigma/2}$}, where $\Delta_N$ is the Neumann Laplacian, defined as the generator of the heat semigroup, see \cite{HM,Functionspace} for example, without involving the boundary from definition. Similarly, one can define another class of Sobolev spaces $H^\sigma_D(K)$ using the Dirichlet Laplacian $\Delta_D$ or the related Brownian motion killed at the boundary.

It is of interest to compare the distinct classes of Sobolev spaces. As one of the main results in this paper, we obtain a full comparison of these spaces, see Theorem \ref{thm32}. Below is the part concerning the relationship between $H^\sigma(K)$ and $H^\sigma_N(K)$.

\begin{thm}
For $k\in\mathbb{Z}_+$, define $\mathcal{H}'_{k-1}=\{f:\Delta^k f=constant,\int_K fd\mu=0\}$, and $d_S$ be the spectral dimension of $\Delta_N$. Then, for $\sigma\in (2k-\frac{d_S}{2},2k+2-\frac{d_S}{2})\cap\mathbb{R}_+$,
	\[H^\sigma(K)=H^{\sigma}_{N}(K)\oplus \mathcal{H}'_{k-1};\]
	for $\sigma=2k+2-\frac{d_S}{2}$, $H^{\sigma}_{N}(K)$ is not closed in $H^\sigma(K)$ and
	\[H^\sigma(K)=cl\big(H^{\sigma}_{N}(K)\big)\oplus \mathcal{H}'_{k-1}.\]
\end{thm}

From Theorem \ref{thm32}, one can easily observe the critical orders $2k+\frac{d_S}{2}$ and $2k+2-\frac{d_S}{2}$ with $k\in \mathbb{Z}_+$, where $d_S$ is the spectral dimension. See \cite{KL} for a discussion on the spectral dimension of $\Delta_N$ or $\Delta_D$. In fact, we have $\frac{d_S}{2}=\frac{d_H}{1+d_H}$, where $d_H$ is the Hausdorff dimension of $K$ with respect to the effective resistance metric(see \cite{book1,book2}), and $d_W:=1+d_H$ is the walk dimension(see \cite{HK,KS}). In particular, on the Sierpinski gasket $\mathcal{SG}$ equipped with the standard energy and measure, we have $d_H=\frac{\log 3}{\log5/3}$, $d_W=\frac{\log5}{\log5/3}$, $d_S=\frac{2\log3}{\log5}$, and the critical orders are $\frac{\log3}{\log5}+2\mathbb{Z}_+$ and $2-\frac{\log3}{\log5}+2\mathbb{Z}_+$.\\

The next topic of this paper is the characterization of the Sobolev spaces, and we are interested in the role of the critical orders. In fact, as a particular case, we usually write a function in $dom\mathcal{E}$ (which equals $H^1_N(K)$, and also $H^1(K)$ by Theorem \ref{thm32}) as a series of tent functions, which is helpful in various cases, for example the trace theorem  of Jonsson \cite{Jonsson2} for $\mathcal{SG}$. We will extend this characterization to general Sobolev spaces. See Theorem \ref{decompositionthm}.

Here we use the Sierpinski gasket $\mathcal{SG}$ as an illustration. We have for $0\leq\sigma<\frac{\log3}{\log5}$, each function $f$ in $H^\sigma(\mathcal{SG})$ admits a unique expansion
\[
f=C+\sum_{w\in W_*}\sum_{i=1}^3 a_{wi}\chi_{_{F_{wi}\mathcal{SG}}},\quad\text{ with } C\in \mathbb{R}, a_{wi}\in \mathbb{R}, \sum_{i=1}^Na_{wi}=0,
\]
and $\|f\|_{H^\sigma(\mathcal{SG})}\asymp \big(|C|^2+\sum_{m=0}^\infty 3^{-m}5^{\sigma m}\sum_{w\in W_m}\sum_{i=1}^3|a_{wi}|^2\big)^{1/2}$;

for  $\frac{\log3}{\log5}<\sigma<2-\frac{\log3}{\log5}$, $f$ adimits \[
f=h+\sum_{x\in V_*\setminus V_0} c_x\psi_x,\quad \text{ with } h\in \mathcal{H}_0, c_x\in\mathbb{R},
\]
and $\|f\|_{H^\sigma(\mathcal{SG})}\asymp \big(\|h\|^2_{L^2(\mathcal{SG})}+\sum_{m=0}^\infty 3^{-m}5^{\sigma m}\sum_{x\in V_{m+1}\setminus V_m}|c_x|^2\big)^{1/2}$; while the expansions can be extended to the case $\sigma\geq 2$ by repeatedly applying the Green's operator. As an application of this result, we invite readers to refer to \cite{tracelower} for the trace spaces of $H^\sigma(\mathcal{SG})$ onto the bottom line segment of $\mathcal{SG}$, which extends Jonsson's work on $dom\mathcal{E}$ \cite{Jonsson2}. {Readers may also compare our decomposition theorem with the theorem of multi-harmonic splines on p.c.f. fractals \cite{spline}.}\\

Lastly, at the mean time of our development, we study various other related characterizations of Sobolev spaces as well.

In 1996, still for the {Sierpinski gasket} $\SG$,  Jonsson\cite{Jonsson1} obtained that the standard energy domain $dom\mathcal{E}$ can be characterized to be a Besov type space $B^{2,\infty}_1(\mathcal{SG})$. Later, the result was generalized to nested fractals by Pietruska-Pa\l uba\cite{Pi} in terms of Euclidean metric and normalized Hausdorff measure, and to p.c.f. fractals by Hu and Wang\cite{HW} by using the effective resistance metric and an associated $d$-regular measure instead.

In 2005, Hu and Z\"{a}hle\cite{HM} studied $H^\sigma_N(K)$, which generalizes $dom\mathcal{E}$, on general metric measure spaces, with the assumption of two-sided heat kernel estimates, which was confirmed to be true for p.c.f. fractals in terms of the effective resistance metric by Hambly and Kumagai, Kumagai and Sturm (see \cite{HK,KS} where probabilistic techniques are much involved). The combination of results in \cite{HM,HK,KS} yields that, for a p.c.f. fractal $K$ with a regular harmonic structure,  $H^\sigma_N(K)$ is equivalent to a Besov type space $B_\sigma^{2,2}(K)$  for $0<\sigma<1$, as well as $H^1_N(K)=dom\mathcal{E}=B_1^{2,\infty}(K)$. We refer to \cite{CG1, CG, GL1} for Besov type characterizations of Sobolev spaces with $0<\sigma<1$ on general metric measure spaces, related to the heat kernel estimates,  and \cite{CG, GKS} for related interpolation results.

It is of interest to find a direct analytical way to characterize a similar Besov type characterization of $H^\sigma(K)$ on p.c.f. fractals. As a particular situation, \cite{Jonsson1,Pi,HW} provided the analytical approach for $dom \mathcal{E}=B^{2,\infty}_1(K)$.

In this paper, we will provide a purely analytic method to show that $H^\sigma(K)=H^\sigma_N(K)=B^{2,2}_\sigma(K)$ for $0<\sigma<1$, without using heat kernel estimate assumption, see Theorem \ref{thm48}. Also, we will take care of the Besov spaces $\Lambda^{2,2}_\sigma(K)$, and the extension $\tilde{\Lambda}^{2,2}_\sigma(K)$ introduced by Strichartz \cite{Functionspace}, and prove that they are all equivalent to $H^\sigma(K)$ for suitable $\sigma$. See Theorem \ref{thm46}, where we show
\[H^\sigma(K)=\Lambda^{2,2}_\sigma(K)\text{ for }\frac{d_S}{2}<\sigma<1\text{,
 and }
H^\sigma(K)=\tilde{\Lambda}^{2,2}_\sigma(K)\text{ for }\frac{d_S}{2}<\sigma<2.\]
We also introduce a new class of Besov type spaces $\Gamma^\sigma(K)$ that equivalent to $H^\sigma(K)$ based on the cell graphs approximating $K$, see Theorem \ref{thm61}. In particular, the case $1<\sigma<2$ are not dealt with in \cite{HM}. 

At the end of the introduction, we briefly show the structure of this paper.

In Section 2, for a p.c.f. fractal $K$, we collect some notations and definitions, including the effective resistance metric $R$, the $d_H$-regular measure $\mu$, and the Sobolev spaces $H^\sigma_D(K)$, $H^\sigma_N(K)$ and $H^\sigma(K)$.

In Section 3, we focus on the comparison theorem of various Sobolev spaces on $K$. We will provide a detailed proof of Theorem \ref{thm32}.

In Section 4, we introduce the main theorems concerning the characterizations of Sobolev spaces, including the atomic decompositions (Theorem \ref{decompositionthm}), and Besov type characterizations (Theorem \ref{thm46} and \ref{thm48}), but postpone the proofs to later sections.

From Section 5 to Section 7, we focus on the characterization of Sobolev spaces with orders $0\leq\sigma<1$. In Section 5, we introduce the Besov type spaces $\Gamma^\sigma(K)$ based on the cell graph energies, and discuss the decomposition in terms of Haar functions. In Section 6, we introduce the notion of so-called smoothed Haar functions, and use it as a key tool to prove the characterizations of Sobolev spaces (Theorem \ref{thm61}, \ref{thm69} and \ref{thm48}). Lastly, in Section 7, we finish the proof of atomic decomposition theorem (Theorem \ref{decompositionthm}) for $0\leq\sigma<1$. At the same time we prove Theorem \ref{thm46}(a) and Theorem \ref{thm59}.

Section 8 is parallel to Section 5 to 7, dealing with the characterization of $H^\sigma(K)$ with higher orders $1\leq\sigma<2$. Since the idea is very similar, we will only provide the key lemmas in this section, and sketch the proof.

Throughout the paper, we always use the notation $f\lesssim g$ if there is a constant $C>0$ such that $f\leq Cg$, and write $f\asymp g$ if $f\lesssim g$ and $g\lesssim f$.

\section{The Dirichlet forms and Sobolev spaces on p.c.f. self-similar sets}
The analysis on p.c.f. self-similar fractals was originally developed by Kigami in \cite{k2,k3}. For convenience of readers, in this section, we will first briefly recall the constructions of Dirichlet forms and Laplacians on p.c.f.  fractals, then introduce the definition of associated Sobolev spaces. Interested readers please refer to  \cite{Hinz,sw} for $p$-energy and corresponding $L^p$ Sobolev spaces on fractals as an extension.

Let $\{F_i\}_{i=1}^N$ be a finite collection of contractions on a complete metric space $(X,d)$. The self-similar set associated with the \textit{iterated function system (i.f.s.)} $\{F_i\}_{i=1}^N$ is the unique compact set $K\subset X$ satisfying
\[K=\bigcup_{i=1}^N F_iK.\]
Each copy $F_iK$ is called a \textit{$1$-cell} of $K$. For $m\geq 1$, we define $W_m=\{1,\cdots,N\}^m$ the collection of \textit{words} of length $m$, and for each $w\in W_m$, denote
\[F_w=F_{w_1}\circ F_{w_2}\circ\cdots \circ F_{w_m}.\]
The set $F_wK$ is called a \textit{$m$-cell} of $K$. Set $W_0=\emptyset$, and let $W_*=\bigcup_{m\geq 0} W_m$ be the collection of all finite words. For $w=w_1w_2\cdots w_m\in W_*\setminus W_0$, we write $w^*=w_1w_2\cdots w_{m-1}$ by deleting the last letter of $w$.

Define the shift space $\Sigma=\{1,2,\cdots,N\}^{\mathbb{N}}$ . There is a continuous surjection $\pi: \Sigma\rightarrow K$ defined by
\[\pi(\omega)=\bigcap_{m\geq 1}F_{[\omega]_m}K,\]
where for $\omega=\omega_1\omega_2\cdots$ in $\Sigma$ we write $[\omega]_m=\omega_1\omega_2\cdots \omega_m\in W_m$ for each $m\geq 1$. Let
\[\mathcal{C}_K=\bigcup_{i\neq j}F_iK\cap F_jK,\quad \mathcal{C}=\pi^{-1}(\mathcal{C}_K),\quad \mathcal{P}=\bigcup_{n\geq 1}\sigma^n \mathcal{C},\]
where $\sigma$ is the shift map define as $\sigma(\omega_1\omega_2\cdots)=\omega_2\omega_3\cdots$. $\mathcal{P}$ is called the \textit{post-critical set}. Call $K$  a \textit{p.c.f. self-similar fractal} if $\#\mathcal{P}<\infty$. In what follows, we always assume that $K$ is a connected p.c.f.  fractal.

Let $V_0=\pi(\mathcal{P})$ and call it the \textit{boundary} of $K$. For $m\geq 1$, we always have $F_w K\cap F_{w'}K\subset F_w V_0\cap F_{w'}V_0$ for any $w\neq w'\in W_m$. Denote $V_m=\bigcup_{w\in W_m}F_wV_0$ and let $l(V_m)=\{f: f \text{ maps } V_m \text{ into } \mathbb{R}\}$. Write $V_*=\bigcup_{m\geq 0}V_m$.

Let $H=(H_{pq})_{p,q\in V_0}$ be a symmetric linear operator(matrix). $H$ is called a \textit{(discrete) Laplacian} on $V_0$ if  $H$ is non-positive definite; $Hu=0$ if and only if $u$ is constant on $V_0$; and $H_{pq}\geq 0$ for any $p\neq q\in V_0$.
Given a Laplacian $H$ on $V_0$ and a vector $\bm{r}=\{r_i\}_{i=1}^N$ with $r_i>0$, $1\leq i\leq N$, define the \textit{(discrete) Dirichlet form} on $V_0$ by
$$\mathcal{E}_0(f,g)=-(f,Hg),$$
and inductively {{on $V_m$ by}}
$$\mathcal{E}_m(f,g)=\sum_{i=1}^Nr^{-1}_i\mathcal{E}_{m-1}(f\circ F_i, g\circ F_i), m\geq 1,$$
for $f, g\in l(V_m)$. Write $\mathcal{E}_m(f,f)=\mathcal{E}_m(f)$ for short.

Say $(H,\bm{r})$ is a \textit{harmonic structure} if for any $f\in l(V_0)$,
\[\mathcal{E}_0(f)=\min\{\mathcal{E}_1(g):g\in l(V_1),g|_{V_{0}}=f\}.\]
In this paper, we will always assume that there exists a harmonic structure associated with $K$, and in addition, $0<r_i<1$ for all $1\leq i\leq N$. Call $(H,\bm{r})$ a \textit{regular harmonic structure} on $K$.

Now for each $f\in C(K)$, the sequence $\{\mathcal{E}_m(f)\}_{m\geq 0}$ is nondecreasing. Let
\begin{equation*}
\mathcal{E}(f,g)=\lim_{m\to\infty} \mathcal{E}_m(f,g) \text{ and }
dom\mathcal{E}=\{f\in C(K):\mathcal{E}(f)<\infty\},
\end{equation*}
where $f,g\in C(K)$ and  we write $\mathcal{E}(f):=\mathcal{E}(f,f)$ for short. Call $\mathcal{E}(f)$ the \textit{energy} of $f$.
It is known that $(\mathcal{E},dom\mathcal{E})$ turns out to be a local regular Dirichlet form on $L^2(K,\mu)$ for any  Randon measure $\mu$ on $K$. 

An important feature of the form $(\mathcal{E},dom\mathcal{E})$ is the \textit{self-similar identity}
\begin{equation}\label{eq21}
\mathcal{E}(f,g)=\sum_{i=1}^Nr_i^{-1}\mathcal{E}(f\circ F_i, g\circ F_i), \quad \forall f,g\in dom\mathcal{E}.
\end{equation}
Furthermore, denote $r_w=r_{w_1}r_{w_2}\cdots r_{w_m}$ for each $w\in W_m, m\geq 0$. Then for  $m\geq 1$, we have
\begin{equation*}
\mathcal{E}_m(f,g)=\sum_{w\in W_m} r_w^{-1}\mathcal{E}_0(f\circ F_w, g\circ F_w),\quad \mathcal{E}(f,g)=\sum_{w\in W_m} r_w^{-1}\mathcal{E}(f\circ F_w, g\circ F_w).
\end{equation*}

See \cite{book1} and \cite{book2} for details and any unexplained notations.

\subsection{Resistance metric and self-similar measure}
To study the Besov spaces on $K$, we need a suitable metric and a comparable measure. Instead of the original $d$, a natural choice of metric is the effective resistance metric $R(\cdot,\cdot)$ \cite{book1}, which matches the form $(\mathcal{E},dom\mcE)$.

\begin{definition}\label{def21}
For $x,y\in K$, the effective resistance metric $R(x,y)$ between $x$ and $y$ is defined by
	\[R(x,y)^{-1}=\min\{\mathcal{E}(f):f\in dom\mcE,f(x)=0,f(y)=1\}.\]
\end{definition}

It is known that $R$ is indeed a metric on $K$ which is topologically equivalent to the metric $d$, and for $w\in W_*$, we always have $diam(F_wK)\asymp r_w$, where $diam(F_wK)=\max\{R(x,y):x,y\in F_wK\}$.

We will always choose {the} following self-similar measure $\mu$ on $K$.

\begin{definition}\label{def22}
	Let $\mu$ be the unique self-similar measure on $K$ satisfying
	\[\mu=\sum_{i=1}^N r_i^{d_H}\mu\circ F_i^{-1},\]
	and $\mu(K)=1$, where $d_H$ is determined by the equation
	$\sum_{i=1}^N r_i^{d_H}=1.$
\end{definition}

For $x\in K$, $\rho>0$, denote $B(x,\rho)=\{y\in K:R(x,y)<\rho\}$  the \textit{ball} centered at $x$ with radius $\rho$. The measure $\mu$ is comparable  with $R$ in the following sense.
\begin{proposition}\label{prop23}
	For any $x\in K$, $0<\rho\leq 1$, we have $\mu\big(B(x,\rho)\big)\asymp \rho^{d_H}$.
\end{proposition}

Before proving the proposition, we  introduce some notation and  lemma, which are of equal importance.

\begin{definition}\label{def24}
	For $0<t\leq 1$,  define $$\Lambda(t)=\{w\in W_*:r_w\leq t<r_{w^*}\}.$$
	In particular, set $r=\min_{i=1}^N r_i$, and let $\Lambda_m=\Lambda(r^m)$ for $m\geq 0$.
\end{definition}

Obviously, $\{\Lambda_m\}_{m\geq 0}$ provides a nested partition of $K$ satisfying that
\[K=\bigcup_{w\in\Lambda_m}F_wK,\quad F_{w}K\cap F_{w'}K\subset F_{w}V_0\cap F_{w'}V_0,\forall w, w'\in \Lambda_m,\]
and thus the self-similar identity (\ref{eq21}) extends to
\begin{equation*}
\mathcal{E}(f,g)=\sum_{w\in \Lambda_m} r_w^{-1}\mathcal{E}(f\circ F_w, g\circ F_w),\quad\forall f,g\in dom\mathcal{E}.
\end{equation*}

\begin{lemma}\label{lemma25}
	For any $x\in K$ and $r^{m+1}<\rho\leq r^{m}$ with $m\geq 0$, we have a uniform bound $M$ over the number of cells intersecting $B(x,\rho)$, i.e.,
	\[\#\{w\in \Lambda_m: F_wK\cap B(x,\rho)\neq \emptyset\}\leq M.\]
\end{lemma}

\textit{Proof.} First, for any $y\in K$, we have a uniform control over the number of cells containing $y$.

\textit{Claim 1: $\#\{w\in \Lambda_m: y\in F_wK\}\leq \#\mathcal{C}$.}

Let $F_{w'}K$ be the smallest cell containing $\bigcup_{w\in\Lambda_m, y\in F_wK}F_wK$, and denote $y'=F_{w'}^{-1}(y)$. Clearly if $\#\{w\in \Lambda_m: y\in F_wK\}\geq 2$, we have $y'\in \mathcal{C}_K$, and
\[\#\{w\in \Lambda_m: y\in F_wK\}\leq\#\pi^{-1}(y')\leq \#\mathcal{C}.\]
Claim 1 is proved.

\textit{{Claim 2:  There exists $C>0$, such that for any $m\geq 1$ and $w,w'\in \Lambda_m$ satisfying $F_wK\cap F_{w'}K=\emptyset$, we have $R(x,y)\geq C r^m$, $\forall x\in F_wK, y\in F_{w'}K$.}}

In fact, for such $w,w',x,y$, choose a function $f\in dom\mathcal{E}$ which takes constant $0$ in $F_wK$, constant $1$ in ${\bigcup\{F_{w''}K:w''\in \Lambda_m,F_{w''}K\cap F_wK=\emptyset\}}$, and is harmonic outside. Then clearly  we have
\[\begin{aligned}
R(x,y)^{-1}\leq \mathcal{E}(f)&= \sum_{\tilde{w}\in \Lambda_m}r_{\tilde{w}}^{-1}\mathcal{E}(f\circ F_{\tilde{w}})\\
&=\sum_{{\tilde{w}}\in \Lambda_m, F_{{\tilde{w}}}K\cap F_{w}K\neq \emptyset}r_{{\tilde{w}}}^{-1}\mathcal{E}(f\circ F_{{\tilde{w}}})\lesssim \#\mathcal{C}\#V_0r^{-m},
\end{aligned}\]
where we use Claim 1 in the last inequality. This gives that $R(x,y)\gtrsim r^m$. Claim 2 follows immediately.

Fix $l$ to be an integer such that $\max\{1,diam(K)\}\cdot 2r^l<C$, where $C$ is the same constant in Claim 2. Assume $m\geq l+1$, then for any $w\in \Lambda_m$ with $F_wK\cap B(x,\rho)\neq\emptyset$, we can pick a point $z\in F_wK\cap B(x,\rho)$ such that $R(x,z)<r^m<\frac 12C r^{m-l}$, and thus for any point $y\in F_wK$, we have $R(x,y) <\frac 12 Cr^{m-l}+r^m\cdot diam(K)<Cr^{m-l}$. By Claim 2, this gives that
\[\begin{aligned}
\bigcup_{w\in \Lambda_m, F_wK\cap B(x,\rho)\neq\emptyset}F_wK\subset  \bigcup_{u\in \Lambda_{m-l}}\{F_uK: \exists &u'\in \Lambda_{m-l}
\text{ such that } \\&x\in F_{u'}K\text{ and } F_uK\cap F_{u'}K\neq \emptyset\}.
\end{aligned}\]
Then the lemma follows from Claim 1 and the fact that for any $u\in \Lambda_{m-l}$, $\#\{w\in \Lambda_m: F_wK\subset F_uK\}$ has a  uniform bound independent of $m$.\hfill$\square$
\vspace{0.2cm}

\textit{Proof of Proposition \ref{prop23}. } Without loss of generality, we assume that $diam(K)=1$.
Let $m\geq 0$ such that $r^{m+1}<\rho\leq r^{m}$. By Lemma \ref{lemma25}, we have
\[\mu\big(B(x,\rho)\big)\leq \sum_{w\in \Lambda_m,F_{w}K\cap B(x,\rho)\neq \emptyset}\mu(K_w)\lesssim Mr^{md_H}\lesssim \rho^{d_H}.\]
On the other hand,  let $w\in \Lambda_{m+1}$ such that $x\in F_wK$. Then since $diam(F_wK)\leq r_w \leq r^{m+1}<\rho$, it is easy to check  $\mu\big(B(x,\rho)\big)\geq \mu(F_wK)=r_w^{d_H}\gtrsim \rho^{d_H}.$ \hfill$\square$
\vspace{0.2cm}

From now on, for simplicity, we will always write $L^2(K)$ instead of $L^{2}(K,d\mu)$, and do similarly for the Sobolev spaces to be defined.

\subsection{Sobolev spaces $H^\sigma(K)$, $H^\sigma_D(K)$ and $H^\sigma_N(K)$} We start from Sobolev spaces $H^\sigma(K)$ with integer orders, then {extend to fractional orders using complex interpolation}. Readers may refer to \cite{pseudo,Jonsson2,distribution,smoothbump,Functionspace} and the references therein for related works, such as bump functions, trace theorems, pseudo-differential operators and a distribution theory on p.c.f. fractals.

We start with the definition of Laplacians.

\begin{definition}\label{def26}
	Let $dom_0\mathcal{E}=\{\varphi\in dom\mathcal{E}:\varphi|_{V_0}=0\}$.

(a). For $f\in dom\mathcal{E}$, say $\Delta f=u$ if
	 \[\mathcal{E}(f,\varphi)=-\int_K u\varphi d\mu , \quad\forall \varphi\in dom_0\mathcal{E}\]

(b). For $f\in dom_0\mathcal{E}$,  say $\Delta_D f=u$ if
\[\mathcal{E}(f,\varphi)=-\int_K u\varphi d\mu, \quad\forall \varphi\in dom_0\mathcal{E},\]

(c). For $f\in dom\mathcal{E}$, say $\Delta_N f=u$ if
\[\mathcal{E}(f,\varphi)=-\int_K u\varphi d\mu, \quad\forall \varphi\in dom\mathcal{E}.\]
\end{definition}

On $L^2(K)$, both $\Delta_D$ and $\Delta_N$ are non-positive definite self-adjoint operators, and $\Delta$ is a closed operator such that $\Delta_D\subset \Delta,\Delta_N\subset \Delta$. In addition, $\Delta_D$ is invertible, {and} $G=-\Delta_D^{-1}$ can be realized with the Green's function $G(x,y)\in C(K\times K)$, i.e.
\[Gf=\int_{K}G(x,y)f(y)d\mu(y).\]
Clearly, we have $-\Delta G f=f,\forall f\in L^2(K)$. For further discussions on the Green's operator $G$, see books \cite{book1} and \cite{book2}.

In the following, we define three different classes of Sobolev spaces, associated with the different Laplacians. One of our main interest in this paper is to clarify their relationships.

\begin{definition}\label{sobolev}
(a). For $k\in \mathbb{Z}_+$, define the Sobolev space $H^{2k}(K)$ as
\[H^{2k}(K)=\{f\in L^2(K): \Delta^j f\in L^2(K) \text{ for all } j\leq k\}\]
with the norm of $f$ given by $\|f\|_{H^{2k}(K)}=\sum_{j=0}^k \|\Delta^j f\|_{L^2(K)}$.

For $0<\theta<1$, $k\in \mathbb{Z}_+$, define $H^{2k+2\theta}(K)$ by using complex interpolation,
\[H^{2k+2\theta}(K)=[H^{2k}(K),H^{2k+2}(K)]_\theta.\]

(b). For $\sigma\geq 0$, define $$H^{\sigma}_{D}(K)=(I-\Delta_D)^{-\sigma/2}L^2(K),$$ with norm $\|f\|_{H^{\sigma}_{D}(K)}=\|(I-\Delta_D)^{\sigma/2}f\|_{L^2(K)}$.

(c).  For $\sigma\geq 0$, define $$H^{\sigma}_{N}(K)=(I-\Delta_N)^{-\sigma/2}L^2(K),$$ with norm $\|f\|_{H^{\sigma}_{N}(K)}=\|(I-\Delta_N)^{\sigma/2}f\|_{L^2(K)}$.
\end{definition}

\section{The critical orders and a comparison theorem}
In this section, we focus on comparing the different Sobolev spaces $H^{\sigma}(K)$, $H^{\sigma}_{D}(K)$ and $H^{\sigma}_{N}(K)$ with $\sigma\geq 0$ defined in Section 2.

Let's start from the {simplest} case, when $\sigma=2k$ for some $k\in \mathbb{Z}_+$.

For $k\in \mathbb{N}$, define the space of \textit{$k$-multiharmonic functions} as
\[\mathcal{H}_{k-1}=\{f:\Delta^k f=0\},\]
and define
\[\mathcal{H}'_{k-1}=\{f:\Delta^k f=constant,\int_K fd\mu=0\}.\]
Both $\mathcal{H}_{k-1}$ and $\mathcal{H}'_{k-1}$ are spaces of dimension $k\#V_0$.
Set $\mathcal{H}_{-1}=\mathcal{H}'_{-1}=\{0\}$ for uniformity.

\begin{proposition}\label{prop31}
Let $k\in \mathbb{Z}_+$. Then we have
\[H^{2k}(K)=H^{2k}_{D}(K)\oplus \mathcal{H}_{k-1},\quad H^{2k}(K)=H^{2k}_{N}(K)\oplus \mathcal{H}'_{k-1}.\]
\end{proposition}
\textit{Proof}. It is easy to verify that $H^{2k}_{D}(K)\cap\mathcal{H}_{k-1} =\emptyset$, $H^{2k}_{N}(K)\cap\mathcal{H}'_{k-1}=\emptyset$, and
$$H^{2k}_{D}(K)\oplus \mathcal{H}_{k-1}\subset H^{2k}(K),\quad H^{2k}_{N}(K)\oplus \mathcal{H}'_{k-1}\subset H^{2k}(K).$$

Next, we show $H^{2k}(K)\subset H^{2k}_{D}(K)\oplus \mathcal{H}_{k-1}$.  Let $f\in H^{2k}(K)$. Define $g=(-1)^kG^kf$ and $h=f-g$. Then, $g\in H^{2k}_{D}(K)$ and $h\in \mathcal{H}_{k-1}(K)$. Thus, $f\in H^{2k}_{D}(K)\oplus \mathcal{H}_{k-1}$.

Last, we show $H^{2k}(K)\subset H^{2k}_N(K)\oplus \mathcal{H}'_{k-1}$. It is well-known that the nonzero eigenvalues of $-\Delta_N$ are bounded away from $0$ and $\ker \Delta_N=constants$, so we may define
\begin{equation}\label{GN}
G_Nu=\sum_{i=1}^\infty \lambda_i^{-1}u_i<u,u_i>,\quad\forall u\in L^2(K),
\end{equation}
where $\{u_i\}_{i=1}^\infty\cup\{1\}$ are eigenfunctions of $-\Delta_N$(chosen to form an orthonormal basis of $L^2(K)$) and $\lambda_i's$ are the corresponding eigenvalues. Let $g=G_N^k\Delta^kf+\int_K fd\mu$ and $h=f-g$, we can see that $\Delta^k h=\int_K\Delta^k fd\mu$ and $\int_Khd\mu=0$, so $h\in \mathcal{H}'_{k-1}$. Thus, $f\in H^{2k}_N(K)\oplus \mathcal{H}'_{k-1}$. \hfill$\square$\vspace{0.2cm}

In this section, {we} will extend Proposition \ref{prop31} to the following Theorem \ref{thm32}. In particular, for the Dirichlet case, we observe a sequence of critical orders that divides $\mathbb{R}_+$ into a sequence of open intervals, such that for $\sigma$ in the $k$-th interval, it holds that $H^{\sigma}(K)=H^{\sigma}_{D}(K)\oplus \mathcal{H}_{k-1}$, and for $\sigma$ being a critical order, the relationship will be more complicated.  The Neumann case is similar, but with a different sequence of critical orders.

\begin{theorem}\label{thm32} Let $d_S=\frac{2d_H}{1+d_H}$ which is the spectral dimension of $K$. We have

(a). For $k\geq 0$ and $\sigma\in (2k-2+\frac{d_S}{2},2k+\frac{d_S}{2})\cap\mathbb{R}_+$,
\[H^\sigma(K)=H^{\sigma}_{D}(K)\oplus \mathcal{H}_{k-1};\]
for $\sigma=2k+\frac{d_S}{2}$, $H^{\sigma}_{D}(K)$ is not closed in $H^\sigma(K)$ and
\[H^\sigma(K)=cl\big(H^{\sigma}_{D}(K)\big)\oplus \mathcal{H}_{k-1}.\]

(b). For $k\geq 0$ and $\sigma\in (2k-\frac{d_S}{2},2k+2-\frac{d_S}{2})\cap\mathbb{R}_+$,
\[H^\sigma(K)=H^{\sigma}_{N}(K)\oplus \mathcal{H}'_{k-1};\]
for $\sigma=2k+2-\frac{d_S}{2}$, $H^{\sigma}_{N}(K)$ is not closed in $H^\sigma(K)$ and
\[H^\sigma(K)=cl\big(H^{\sigma}_{N}(K)\big)\oplus \mathcal{H}'_{k-1}.\]
\end{theorem}

As an immediate consequence of Theorem \ref{thm32}, we have the following useful corollary.
\begin{corollary}\label{coro33}
	$H^1(K)=H^{1}_{D}(K)\oplus\mathcal{H}_0=H^{1}_{N}(K)=dom\mathcal{E}$.
\end{corollary}

We will prove Theorem \ref{thm32} in the rest of this section, and break the proof into several lemmas. Since the proof of part (a) and part (b) are very similar, we will focus on the proof of (a), and sketch the proof of (b) at the same time.

\subsection{The sequence spaces}
As an important tool, we introduce the following sequence spaces. Throughout this section, we always use the symbol $\alpha=(\alpha_1,\alpha_2,\cdots)$ to denote a sequence. Recall that $d_W=1+d_H$ and $d_S=\frac{2d_H}{d_W}$. For $\sigma\geq 0$, we denote by $\lambda=\lambda(\sigma)=r^{(d_H-\sigma d_W)/2}$ for short, where $r=\min_{i=1}^Nr_i$. Let
\[\begin{aligned}
{S}^\sigma=&\big\{\alpha:\{\lambda^m(\alpha_{m+1}-\alpha_m)\}_{m\geq 1}\in l^2\big\},\\
\tilde{{S}}^{\sigma}=&\big\{\alpha:\{\lambda^m\alpha_m\}_{m\geq 1}\in l^2\big\},
\end{aligned}\]
with norms \[\begin{aligned}
\|\alpha\|_{{S}^\sigma}=&|\alpha_1|+\|\lambda^m(\alpha_{m+1}-\alpha_m)\|_{l^2},\\
\|\alpha\|_{\tilde{{S}}^{\sigma}}=&\|\lambda^m\alpha_{m}\|_{l^2},
\end{aligned}\]
respectively.

\begin{lemma}\label{lemma34}
(a). For $\sigma<\frac{d_S}{2}$, we have $S^\sigma=\tilde{S}^{\sigma}$ with $\|\alpha\|_{S^\sigma}\asymp \|\alpha\|_{\tilde{S}^\sigma}$.

(b) For $\sigma>\frac{d_S}{2}$, we have $S^\sigma=constants\oplus \tilde{S}^{\sigma}$ with $\|\alpha\|_{S^\sigma}\asymp|c|+ \|\tilde{\alpha}\|_{\tilde{S}^{\sigma}}$, where $c=\lim_{m\to\infty}\alpha_m$ and $\tilde{\alpha}_m=\alpha_m-c$.

(c). For $\sigma=\frac{d_S}{2}$, $\tilde{S}^\sigma$ is not a closed subspace of $S^\sigma$, and $\tilde{S}^\sigma$ is dense in $S^\sigma$.
\end{lemma}

\textit{Proof.} (a). Let $\sigma<\frac{d_S}{2}$ and so $\lambda<1$. For any $\alpha\in S^\sigma$, by Minkowski inequality, we have
\[\begin{aligned}
\|\alpha\|_{\tilde{S}^\sigma}=\|\lambda^m\alpha_m\|_{l^2}&{= \big\|\lambda^m\alpha_1+\lambda^m\sum_{k=1}^{m-1}(\alpha_{k+1}-\alpha_k)\big\|_{l^2}}\\
&{\lesssim}|\alpha_1|+\big\|\sum_{k=1}^{m-1}\lambda^k\cdot\lambda^{m-k}(\alpha_{m-k+1}-\alpha_{m-k})\big\|_{l^2}\\
&\leq |\alpha_1|+\sum_{k=1}^\infty \lambda^{k}\|\lambda^{m}(\alpha_{m+1}-\alpha_{m})\|_{l^2}\lesssim \|\alpha\|_{S^\sigma}.
\end{aligned}\]
Conversely, for any $\alpha\in \tilde{S}^\sigma$, it is trivial that $\|\alpha\|_{S^\sigma}\lesssim \|\alpha\|_{\tilde{S}^\sigma}$.

(b). Let $\sigma>\frac{d_S}{2}$ and so $\lambda>1$. For any $\alpha\in {S}^\sigma$, it is easy to see the limit $c:=\lim_{m\to\infty} \alpha_m=\alpha_1+\sum_{m=1}^{\infty}(\alpha_{m+1}-\alpha_m)$ exists. Let $\tilde{\alpha}_m=\alpha_m-c$. Then by Minkowski inequality, we have
\[\begin{aligned}
\|\tilde{\alpha}\|_{\tilde{S}^\sigma}=\|\lambda^m\tilde{\alpha}_m\|_{l^2}&=\big\|\lambda^m\sum_{k=m}^{\infty}(\alpha_{k+1}-\alpha_k)\big\|_{l^2}\\
&=\big\|\sum_{k=0}^{\infty}\lambda^{-k}\cdot\lambda^{k+m}(\alpha_{k+m+1}-\alpha_{k+m})\big\|_{l^2}\\
&\leq \sum_{k=0}^\infty \lambda^{-k}\|\lambda^{m}(\alpha_{m+1}-\alpha_{m})\|_{l^2}\lesssim \|\alpha\|_{{S}^{\sigma}}.
\end{aligned}
\]
In addition, it is trivial that $|c|\lesssim \|\alpha\|_{{S}^{\sigma}}$. Thus we have $|c|+\|\tilde{\alpha}\|_{\tilde{S}^\sigma}\lesssim \|\alpha\|_{S^\sigma}$. The other direction of the estimate is trivial.

(c). Let $\sigma=\frac{d_S}{2}$ and so $\lambda=1$. Then $\tilde{S}^\sigma=l^2$ and $S^\sigma=\big\{\alpha:\{(\alpha_{m+1}-\alpha_m)\}_{m\geq 1}\in l^2\big\}$. The claim is obvious.\hfill$\square$\vspace{0.2cm}

Before proceeding to the proof of Theorem \ref{thm32}, we introduce some more notations here. We will write $\bm{\alpha}=\{\alpha^{p}\}_{p\in V_0}$ with each $\alpha^p=(\alpha^p_1,\alpha^p_2,\cdots)$ being a sequence. In addition, let ${\mathcal{S}}^\sigma=({S}^\sigma)^{V_0}$ and $\mathcal{\tilde{S}}^{\sigma}=(\tilde{S}^{\sigma})^{V_0}$, with norms
\[\|\bm{\alpha}\|_{{\mathcal{S}}^\sigma}=\sum_{p\in V_0}\|\alpha^{p}\|_{{S}^\sigma},\quad\|\bm{\alpha}\|_{\mathcal{\tilde{S}}^{\sigma}}=\sum_{p\in V_0} \|\alpha^{p}\|_{\tilde{S}^{\sigma}}.\]

\subsection{Embedding the sequence spaces}
In this subsection, we will embed the sequence spaces ${\mathcal{S}}^\sigma$ or $\tilde{{\mathcal{S}}}^\sigma$ into the Sobolev spaces $H^\sigma(K)$, $H_D^\sigma(K)$ or $H_N^\sigma(K)$. In particular, we will introduce a ``restriction'' map and an ``extension'' map for different cases.

First, we introduce some notations. For $m\geq 0$, let $\Lambda_m=\Lambda(r^m)$ as introduced in Definition \ref{def24} with $r=\min_{i=1}^N r_i$. For each $p\in V_0$, $m\geq 0$, denote
$$\Lambda_{p,m}=\{w\in \Lambda_m: p\in F_wK\}\text{, and }U_{p,m}=\bigcup_{w\in \Lambda_{p,m}}F_wK.$$
Then $\{U_{p,m}\}_{m\geq 0}$ is a decreasing sequence of neighbourhoods of $\{p\}$. Without loss of generality, we assume

 1. $\#\Lambda_{p,1}=\#\pi^{-1}(p)$ for any $p\in V_0$;

 2.  $U_{p,1}\cap U_{q,1}=\emptyset$ for any $p,q\in V_0$;

 3. $\overline{(U_{p,m})^c}\cap U_{p,m+1}=\emptyset$ for any $p\in V_0$ and $m\geq 1$, \\
  \noindent otherwise we replace $r$ by a sufficiently small number.

\begin{lemma}\label{lemma35}
	For each $p\in V_0$, there exist two functions $\phi_{p}$ and $\psi_{p}$ in $\mathcal{H}_0$ such that for any $h\in \mathcal{H}_0$, we have
	\[\int_{K}\phi_p hd\mu=h(p)\text{, and }\int_{K}\psi_p hd\mu=\partial_nh(p).\]
\end{lemma}
\textit{Proof.} Observing that $h\to h(p)$ and $h\to \partial_n h(p)$ are functionals on $\mathcal{H}_0$, the lemma follows immediately from Riesz representation theorem.\hfill$\square$

\begin{definition}\label{def36}
	For $p\in V_0$, $m\geq 1$, let $\phi_{p,m}$ and  $\psi_{p,m}$ be two functions supported in $U_{p,m}$ such that $$\phi_{p,m}=\sum_{w\in \Lambda_{p,m}}r_w^{-d_H}\phi_{F_w^{-1}p}\circ F_w^{-1}\text{ and }\psi_{p,m}=\sum_{w\in \Lambda_{p,m}}r_w^{-1-d_H}\psi_{F_w^{-1}p}\circ F_w^{-1}.$$
	Furthermore, for $f\in L^2(K)$, define
	
	(a). $R^{p}_vf=\big\{\int_K \phi_{p,m}fd\mu\}_{m\geq 1}$ and write $R_vf=\{R_v^{p}f\}_{p\in V_0}$;
	
	(b). $R^{p}_nf=\big\{\int_K \psi_{p,m}fd\mu\big\}_{m\geq 1}$ and write $R_nf=\{R_n^{p}f\}_{p\in V_0}$.
\end{definition}

The operators $R_v$ and $R_n$ will play the role of the ``restriction'' map.

\begin{lemma}\label{lemma37} Let $0\leq\sigma\leq 2$.
	
	(a). The map $R_v$ is bounded from $H^\sigma(K)$ to $\mathcal{S}^\sigma$, and also from $H^\sigma_D(K)$ to $\mathcal{\tilde{S}}^{\sigma}$.
	
	(b). The map $R_n$ is bounded from $H^\sigma(K)$ to ${\mathcal{S}}^{\sigma-{{2}}/{d_W}}$, and also from $H^\sigma_N(K)$ to $\mathcal{\tilde{S}}^{\sigma-{{2}}/{d_W}}$.
\end{lemma}

\textit{Proof.} We only prove (a), and the proof of (b) is essentially the same. First we show $R_v: H^\sigma(K)\rightarrow \mathcal{{S}}^{\sigma}$  is bounded. By complex interpolation, we only need to show it for $\sigma=0$ and $\sigma=2$.

For $\sigma=0$, it follows from the estimates that for any $f\in L^2(K)$ and any $p\in V_0$, we have
\begin{equation}\label{eqn81}\begin{aligned}
\|R_v^{p}f\|_{{S}^0}&\asymp\|R_v^{p}f\|_{\tilde{S}^0}
=\big\|{r^{md_H/2}}\int_{U_{p,m}}\phi_{p,m} f d\mu\big\|_{l^2}
\lesssim \big\|r^{-md_H/2}\int_{U_{p,m}}|f|d\mu\big\|_{l^2}\\
&\lesssim \big\|r^{-md_H/2}\sum_{k=m}^{\infty}r^{kd_H/2}\|f\|_{L^2(U_{p,k}\setminus U_{p,k+1})}\big\|_{l^2}\\
&=\big\|r^{-md_H/2}\sum_{k=0}^{\infty}r^{(m+k)d_H/2}\|f\|_{L^2(U_{p,m+k}\setminus U_{p,m+k+1})}\big\|_{l^2}\\
&=\sum_{k=0}^{\infty}r^{kd_H/2}\big\|\|f\|_{L^2(U_{p,m+k}\setminus U_{p,m+k+1})}\big\|_{l^2}\lesssim \|f\|_{L^2(K)}.
\end{aligned}
\end{equation}
where we use Lemma \ref{lemma34} (a) in the first estimate, and use Cauchy-Schwartz inequality, Minkowski inequality in the remaining estimates.

For $\sigma=2$, let $\tilde{\phi}_{p,m}=\phi_{p,m+1}-\phi_{p,m}, \forall p\in V_0, m\geq 1$. For each $w\in \Lambda_{p,m}$, immediately we have $\int_{F_wK} \tilde{\phi}_{p,m}hd\mu=0$, for each $h$ harmonic in $F_wK$. As a result, by using Gauss-Green's formula on each $F_wK$, for any $f\in H^2(K)$, we have
\[\int_{F_wK} \tilde{\phi}_{p,m}fd\mu=\int_{F_wK} \tilde{\phi}_{p,m}(f-h)d\mu=-\int_{F_wK} G_w\tilde{\phi}_{p,m}\cdot \Delta fd\mu,\]
where $h$ is {harmonic in $F_wK$ with $h|_{F_wV_0}=f|_{F_wV_0}$}, and $G_w$ is the local Green's function associated with $F_wK$. Define $$\phi'_{p,m}=-\sum_{w\in \Lambda_{p,m}}G_w\tilde{\phi}_{p,m}.$$
Then, it is easy to see that
\[(R_v^pf)_{m+1}-(R_v^pf)_{m}=\int_{U_{p,m}} \phi'_{p,m}\cdot \Delta fd\mu.\]
In addition, we have the estimate $\|\phi'_{p,m}\|_{L^\infty(U_{p,m})}\lesssim r^m$. The result for $\sigma=2$ then follows a by similar estimate as (\ref{eqn81}).

For the boundedness of $R_v:H^\sigma_D(K)\to \mathcal{\tilde{S}}^{\sigma}$, we still use the complex interpolation for $0\leq \sigma\leq 2$. For $\sigma=0$, it follows  immediately since $H^0_D(K)=L^2(K)$. For $\sigma=2$, we only need to notice that for any $p\in V_0$, $f\in H^2_D(K)$, we always have $\lim_{m\rightarrow\infty}r^{-md_H/2}\int_{U_{p,m}}|f|d\mu=0$, so that $R_vf\in \mathcal{\tilde{S}}^2$ and $\|R_vf\|_{\mathcal{S}^2}\asymp \|R_vf\|_{\mathcal{\tilde{S}}^2}$ by Lemma \ref{lemma34} (b). \hfill$\square$\\

\noindent\textbf{Remark.} By a routine discussion, one can easily see $\lim\limits_{m\to\infty} (R^p_vf)_m=\#\pi^{-1}(p)f(p)$ and $\lim\limits_{m\to\infty} (R^p_nf)_m=\partial_n f(p)$ for any $f\in H^2(K)$ and $p\in V_0$. \\

In the next lemma, we construct maps $E_v,E_n$ that play the role of ``extension'' map.

\begin{lemma}\label{lemma38} Let $0\leq \sigma\leq 2$.
	
	(a). There exists a bounded map $E_v:\mathcal{S}^\sigma\to H^\sigma(K)$ satisfying $R_v\circ E_v=id$.
	
	(b). There exists a bounded map $E_n:\mathcal{S}^{\sigma-{{2}}/{d_W}}\to H^\sigma(K)$ satisfying $R_n\circ E_n=id$.
\end{lemma}

\textit{Proof.} (a). For $p\in V_0$ and $m\geq 1$, we choose a function $g_{p,m}\in H^2(K)$ such that
\[
g_{p,m}|_{U_{p,m+1}}=\big(\#\pi^{-1}(p)\big)^{-1},\quad  g_{p,m}|_{K\setminus U_{p,m}}=0,
\]
and
\[(R^{p}_vg_{p,m})_k=\begin{cases}
0,\text{ for }k\leq m,\\
1,\text{ for }k\geq m+1.
\end{cases}\]
In addition,  for each $p\in V_0$, we can guarantee
\[\sup_{m\geq 1} r^{-m(d_H/2+1)}\|\Delta g_{p,m}\|_{L^2(K)}<\infty, \quad \sup_{m\geq 1}\|g_{p,m}\|_{L^\infty(K)}<\infty.\]
In fact, we can construct proper $g_{p,m}|_{F_wK\setminus U_{p,m+1}}$ for each $w\in \Lambda_{p,m}$ respectively, and only finitely many different cases of $F_wK\setminus U_{p,m+1}$ ({up to some contraction mapping $F_v$}) will occur. The estimate for $\|g_{p,m}\|_{L^\infty(K)}$ is immediate, and the estimate for $\|\Delta g_{p,m}\|_{L^2(K)}$ follows from the scaling property of $\Delta$.

For each $\bm{\alpha}\in \mathcal{S}^\sigma$, define
\[E_v\bm{\alpha}=h+\sum_{p\in V_0}\sum_{m=1}^\infty (\alpha^{p}_{m+1}-\alpha^{p}_{m})g_{p,m},\]
where $h$ is a harmonic function with boundary values $h(p)=\big(\#\pi^{-1}(p)\big)^{-1}\alpha_1^{p}$, $\forall p\in V_0$.

We will show that $E_v:\mathcal{S}^\sigma\to L^2_\sigma(K)$ is bounded. By complex interpolation, it is enough to show it for $\sigma=0$ and $\sigma=2$.  For $\sigma=0$, by using Minkowski inequality, for each $\alpha\in \mathcal{S}^0$, we have the estimate
\[\begin{aligned}
\|E_v\bm{\alpha}\|_{L^2(K)}&\leq \sum_{p\in V_0}\big\|\sum_{k=1}^\infty (\alpha^{p}_{k+1}-\alpha^{p}_k)g_{p,k}\big\|_{L^2(U_{p,1})}+\|\bm{\alpha}\|_{\mathcal{S}^0}\\
&=\sum_{p\in V_0}\big\|\|\sum_{k=1}^m (\alpha^{p}_{k+1}-\alpha^{p}_{k})g_{p,k}\|_{L^2(U_{p,m}\setminus U_{p,m+1})}\big\|_{l^2}+\|\bm{\alpha}\|_{\mathcal{S}^0}\\
&\lesssim \sum_{p\in V_0}\big\| r^{md_H/2}\sum_{k=1}^{m}|\alpha^{p}_{k+1}-\alpha^{p}_{k}|\big\|_{l^2}+\|\bm{\alpha}\|_{\mathcal{S}^0}\\
&=\sum_{p\in V_0}\big\|\sum_{k=0}^{m-1} r^{kd_H/2}\cdot r^{(m-k)d_H/2} |\alpha^{p}_{m-k+1}-\alpha^{p}_{m-k}|\big\|_{l^2}+\|\bm{\alpha}\|_{\mathcal{S}^0}\\&\lesssim \|\bm{\alpha}\|_{\mathcal{S}^0}.
\end{aligned}\]
For $\sigma=2$, the proof is easy, noticing each $\Delta g_{p,m}$ is locally supported on $U_{p,m}\setminus U_{p,m+1}$.

Lastly, it is direct to check that $R_v(E_v\bm{\alpha})=\bm{\alpha}$ on $\mathcal{S}^0$.

The proof of (b) is essentially the same. The main difference is that we construct $g_{p,m}$ to be harmonic in each cell of $U_{p,m+1}$ with desired normal derivative at $p$. We omit the details.\hfill$\square$\\

As an immediate consequence of Lemma \ref{lemma38}, and by the remark after Lemma \ref{lemma37}, we have the following lemma.

\begin{lemma}\label{lemma39} Let $0\leq \sigma\leq 2$, and $E_v, E_n$ be the map defined in Lemma \ref{lemma38}. Then

	(a). $E_v$ is bounded from $\mathcal{\tilde{S}}^\sigma\to H^\sigma_D(K)$.
	(b). $E_n$ is bounded from $\mathcal{\tilde{S}}^{\sigma-2/d_W}$ to $H^\sigma_N(K)$.
\end{lemma}

\subsection{Proof of the comparison theorem} In this part, we come to the proof of Theorem \ref{thm32}.
We need to use the following simple fact, which can be easily derived from the property of interpolation functors.
\begin{lemma}\label{lemma310}
	Let $(Z_1,Z_2)$ be an interpolation couple with $Z_1=X_1\oplus Y_1$, $Z_2=X_2\oplus Y_2$, and  $(X_1+X_2)\cap (Y_1+Y_2)=\{0\}$. Then we have
	$$[Z_1,Z_2]_\theta=[X_1,X_2]_\theta\oplus [Y_1,Y_2]_\theta,\quad\forall 0<\theta<1.$$
\end{lemma}

The following lemma  concludes what we have got in the last two subsections.

\begin{lemma}\label{lemma311} Let $0\leq\sigma\leq 2$. Define
	\[\ker_\sigma R_v=\{f\in H^\sigma(K):R_vf=0\}\text{ and }\ker_\sigma R_n=\{f\in H^\sigma(K):R_nf=0\}.\]
	Then we have

(a). $H^\sigma(K)=E_v\mathcal{S}^\sigma\oplus \ker_\sigma R_v$ and $H^\sigma_D(K)=E_v\mathcal{\tilde{S}}^{\sigma}\oplus \ker_\sigma R_v$.
	
(b). $H^\sigma(K)=E_n\mathcal{S}^{\sigma-{2}/{d_W}}\oplus \ker_\sigma R_n$ and $H^\sigma_N(K)=E_n\mathcal{\tilde{S}}^{\sigma-{2}/{d_W}}\oplus \ker_\sigma R_n$.
\end{lemma}

\textit{Proof.} (a). The first identity is obvious
by Lemma \ref{lemma37} and Lemma \ref{lemma38}. In addition, we can see that \begin{equation}\label{eqn32}
\ker_\sigma R_v=[\ker_0 R_v,\ker_2 R_v]_{\sigma/2},\quad\forall 0<\sigma<2,
\end{equation}
by applying Lemma \ref{lemma310}.

We can also see from the first identity that $\ker_\sigma R_v=\{f\in H^\sigma_D(K):R_vf=0\}$ for $\sigma=0,2$ using the remark after Lemma \ref{lemma37}. Applying Lemma \ref{lemma37} and Lemma \ref{lemma39}, we then have the second identity holds for $\sigma=0,2$. Since $H^\sigma_D(K)$ and $\mathcal{\tilde{S}}^\sigma$ are stable under complex interpolation, we have
\[H^\sigma_D(K)=E_v\mathcal{\tilde{S}}^\sigma\oplus[\ker_0 R_v,\ker_2 R_v]_{\sigma/2},\]
by applying Lemma \ref{lemma310}. The second identify then follows from (\ref{eqn32}) immediately.

The proof of (b) is the same.\hfill$\square$\vspace{0.2cm}

\textit{Proof of Theorem \ref{thm32}}. (a). For $0\leq \sigma\leq 2$, the claims are easy {consequences} of Lemma \ref{lemma311}. For $0\leq\sigma<\frac{d_S}{2}$, the result follows from Lemma \ref{lemma311} and Lemma \ref{lemma34} (a). For $\sigma=\frac{d_S}{2}$, the result follows from Lemma \ref{lemma311} and Lemma \ref{lemma34} (c). For $\frac{d_S}{2}<\sigma\leq 2$, one have $\mathcal{H}_0\cap H^\sigma_D(K)=\{0\}$ since obviously $R_v\mathcal{H}_0\cap R_vH^\sigma_D(K)=\{0\}$ and $\ker_\sigma R_v\cap \mathcal{H}_0=\{0\}$. The result follows from Lemma \ref{lemma311}, Lemma \ref{lemma34} (b) and codimension counting.

The result for $2k\leq\sigma\leq 2k+2$ follows from the fact that
\begin{equation}\label{eqn34}
H^\sigma_D(K)=G^kH^{\sigma-2k}_D(K) \text{, and }H^\sigma(K)=\mathcal{H}_{k-1}\oplus G^kH^{\sigma-2k}(K).
\end{equation}
To see the second equality, we follow a similar proof as Proposition \ref{prop31} for $\sigma=2k,2k+2$ and then apply Lemma \ref{lemma310}.

(b). For $0\leq\sigma\leq 2$, the proof is the same as (a). For $2k\leq\sigma\leq 2k+2$, similar to (a), we have
\[H^\sigma_N(K)=constants\oplus G_N^kH^{\sigma-2k}_N(K)\text{, and }H^\sigma(K)=\mathcal{H}'_{k-1}\oplus \big(constants\oplus G_N^kH^{\sigma-2k}(K)\big).\]
\hfill$\square$

\section{The atomic decomposition and other Besov type characterizations}
It is of interest to see what kind of role the critical orders in Theorem \ref{thm32} play. It is well-known that $\frac{d_S}{2}$ is a critical order of continuity of functions in Sobolev spaces $H^\sigma(K)$. Also, we expect that $2-\frac{d_S}{2}$ is a critical order concerning H\"{o}lder continuity of functions in $H^\sigma(K)$. The forthcoming atomic decomposition theorem will provide a more superising explanation.

\subsection{The atomic decomposition}

We denote $\chi_A$ for the characteristic function on a set $A$ contained in $K$, i.e. $\chi_A(x)=1$ if $x\in A$ and $\chi_A(x)=0$ if $x\notin A$.

We denote $\psi_x$ for the tent function at a point $x\in V_*\setminus V_0$. To be more precise, for $m\in\mathbb{N}$ and $x\in V_m\setminus V_{m-1}$, we define $\psi_x(y)=1$ if $y=x$, $\psi_x(y)=0$ if $y\in V_{m}\setminus \{x\}$, and extend $\psi_x$ to be  harmonic  in each cell $F_wK$ with $w\in W_m$.

\begin{theorem}\label{decompositionthm}
	(a). For $k\geq 0$ and $\sigma\in (2k-\frac{d_S}{2},2k+\frac{d_S}{2})\cap \mathbb{R}_+$,  the series
	\begin{equation}\label{eq41}
	f=h+G^kC+\sum_{w\in W_*}\sum_{i=1}^N a_{wi}G^k\chi_{_{F_{wi}K}},
	\end{equation}
	with $h\in\mathcal{H}_{k-1}$, $C\in \mathbb{R}$, $a_{wi}\in \mathbb{R}$, $\sum_{i=1}^Nr_i^{d_H}a_{wi}=0$, is in $H^\sigma(K)$ if and only if
	\[
	\sum_{w\in W_*}r_w^{d_H-(\sigma-2k) d_W} \sum_{i=1}^N|a_{wi}|^2<\infty,
	\]
	with $\|f\|_{H^\sigma(K)}\asymp \big(\|h\|^2_{L^2(K)}+|C|^2+\sum_{w\in W_*}r_w^{d_H-(\sigma-2k) d_W} \sum_{i=1}^N|a_{wi}|^2\big)^{1/2}$.
	In addition, each $f$ in $H^\sigma(K)$ admits a unique expansion of this form.
	
	(b). For $k\geq 0$ and  $\sigma\in(2k+\frac{d_S}{2},2k+2-\frac{d_S}{2})\cap\mathbb{R}_+$, the series \begin{equation}\label{eq42}
	f=h+\sum_{x\in V_*\setminus V_0} c_xG^k\psi_x,
	\end{equation} with $h\in \mathcal{H}_k$ and $c_x\in\mathbb{R}$, is in  $H^\sigma(K)$ if and only if
	\[
	\sum_{w\in W_*}r_w^{d_H-(\sigma-2k) d_W}\sum_{x\in V_1\setminus V_0}|c_{F_wx}|^2<\infty,
	\]
	with $\|f\|_{H^\sigma(K)}\asymp \big(\|h\|^2_{L^2(K)}+\sum_{w\in W_*}r_w^{d_H-(\sigma-2k) d_W}\sum_{x\in V_1\setminus V_0}|c_{F_wx}|^2\big)^{1/2}$. In addition, each $f$ in $H^\sigma(K)$ admits a unique expansion of this form.
\end{theorem}

\noindent\textbf{Remark 1.} Comparing Theorem \ref{decompositionthm} with Theorem \ref{thm32},  we can immediately  get a similar atomic decomposition theorem for $H^{\sigma}_{D}(K)$ with $\sigma\geq 0$. The statement is almost the same with the multiharmonic function term $h$ removed from the expansions in (\ref{eq41}) and (\ref{eq42}).

\noindent\textbf{Remark 2.} There are some other types of atomic decomposition theorems. See Theorem \ref{thm69} and Theorem \ref{thm84} in later sections.
\vspace{0.15cm}

We will prove Theorem \ref{decompositionthm} in the remaining part of this
paper. Due to formula (\ref{eqn34}), we only need to consider orders $0\leq\sigma<2$. We will split the proof into two main parts:  $0\leq \sigma< 1$ and $1\leq \sigma< 2$. The case $1\leq\sigma< 2$ is more or less similar to the case $0\leq\sigma< 1$,  and so  some details will be omitted in the second part.

\subsection{Besov type characterizations}
At the mean time, we will provide some other types of characterizations of $H^\sigma(K)$. Again, we will postpone the proof in later sections. Throughout this paper, we write $\lambda=\lambda(\sigma)=r^{(d_H-\sigma d_W)/2}$ as in Section 3.

We follow Strichartz \cite{Functionspace} to define the following Besov type spaces $\Lambda^{2,2}_\sigma(K)$ and $\tilde{\Lambda}^{2,2}_\sigma(K)$, based on discrete differences. The spaces are defined with different differences for $\frac{d_S}{2}<\sigma<1$ and $\frac{d_S}{2}<\sigma<2$ separately.

For $\frac{d_S}{2}<\sigma<1$, we consider the following spaces.
\begin{definition}\label{def42}
	For $m\geq 0$, we write $V_{\Lambda_m}=\bigcup_{w\in \Lambda_m}F_wV_0$ for short.
	
	(a). For $m\geq 0$, call the graph $G_{v,m}=(V_{\Lambda_m}, E_{v,m})$ with vertice set $V_{\Lambda_m}$ and edge set $E_{v,m}$ a level-$m$ vertex graph, where $\{x,y\}\in E_{v,m}$ if and only if there exists $w\in \Lambda_m$ and $p,q\in V_0$ such that $x=F_wp,y=F_wq$.
	
	(b). For $m\geq 0$, define the difference operator $\nabla_m: C(K)\to l(E_{v,m})$ as
	$$\nabla_mf(\{x,y\})=f(x)-f(y), \quad \forall f\in C(K) \text{ and } \{x,y\}\in E_{v,m}.$$ Notice that we fix one direction for each edge $\{x,y\}$.
\end{definition}

\begin{definition}\label{def43}
	For $\sigma>\frac{d_S}{2}$, define
	\[\Lambda^{2,2}_\sigma(K)=\{f\in C(K):\sum_{m=0}^\infty \lambda^{2m}\|\nabla_mf\|^2_{l^2(E_{v,m})}<\infty\},\]
	with norm $\|f\|_{\Lambda^{2,2}_\sigma(K)}=\big(\|f\|^2_{L^2(K)}+\sum_{m=0}^\infty \lambda^{2m}\|\nabla_mf\|^2_{l^2(E_{v,m})}\big)^{1/2}$.
\end{definition}

\noindent{\textbf{Remark.} We point out that $\Lambda_{\sigma}^{2,2}(K)=constants$ when $\sigma\geq 1$, since then $$\mathcal{E}(f)\lesssim\lim_{m\to\infty} r^{-m}\|\nabla_{m}f\|_{l^2(E_{v,m})}^2=0,\quad\forall f\in \Lambda_{\sigma}^{2,2}(K).$$}

For $\frac{d_S}{2}<\sigma<2$, we consider the following spaces.

\begin{definition}\label{def44}
	 For $m\geq 0$, denote the (discrete) Dirichlet form on $V_{\Lambda_m}$ by
	\[\mathcal{E}_{\Lambda_m}(f,g)=\sum_{w\in\Lambda_m} r_w^{-1}\mathcal{E}_0(f\circ F_w, g\circ F_w), \quad \forall f,g\in l(V_{\Lambda_m}),\]
	and write its corresponding graph Laplacian as $H_{\Lambda_m}:l(V_{\Lambda_m})\to l(V_{\Lambda_m})$, i.e.,
	\[\mathcal{E}_{\Lambda_m}(f,g)=-(f, H_{\Lambda_m}g), \quad\forall f,g\in l(V_{\Lambda_m}).\]
\end{definition}

\begin{definition}\label{def45}
	For $\sigma>\frac{d_S}{2}$, define
	\[\tilde{\Lambda}^{2,2}_\sigma(K)=\{f\in C(K):\sum_{m=1}^\infty r^{2m}\lambda^{2m}\|H_{\Lambda_m}f\|^2_{l^2(V_{\Lambda_m}\setminus V_0)}<\infty\},\]
	with norm $\|f\|_{\tilde{\Lambda}^{2,2}_\sigma(K)}=\big(\|f\|^2_{L^2(K)}+\sum_{m=1}^\infty r^{2m}\lambda^{2m}\|H_{\Lambda_m}f\|^2_{l^2(V_{\Lambda_m}\setminus V_0)}\big)^{1/2}$.
\end{definition}

We will prove the following theorem.

\begin{theorem}\label{thm46}
	(a). For $\frac{d_S}{2}<\sigma<1$, we have $H^\sigma(K)=\Lambda^{2,2}_\sigma(K)$ with $\|\cdot\|_{H^\sigma(K)}\asymp\|\cdot\|_{\Lambda_\sigma^{2,2}(K)}$.
	
	(b). For $\frac{d_S}{2}<\sigma<2$, we have $H^\sigma(K)=\tilde{\Lambda}^{2,2}_\sigma(K)$ with $\|\cdot\|_{H^\sigma(K)}\asymp\|\cdot\|_{\tilde{\Lambda}_\sigma^{2,2}(K)}$.
\end{theorem}

Another class of Besov type spaces on fractals, denoted by $B^{2,2}_\sigma(K)$, given by the double integration, has been widely studied in connection with heat kernel and Dirichlet forms, see {\cite{GL,HiK,HW,HM,Jonsson1,Pi}}. Also see \cite{BV1,BV2} for analogous spaces with slight differences. It was shown in \cite{HM} on general metric measure spaces that $H^\sigma_N(K)$ are identical with $B^{2,2}_{\sigma}(K)$ when $0<\sigma<1$ under the assumption of nice heat kernel estimates. As an application of our results, we will give a {purely} analytical way to prove this result without using the heat kernel estimates.

\begin{definition}\label{def47}
	For $\sigma>0$ and $f\in L^2(K)$, denote
	 \[[f]_{B^{2,2}_\sigma(K)}=\big(\int_{0}^1\frac{dt}{t}\int_{K}\int_{B(x,t)}t^{-d_H-\sigma d_W}|f(x)-f(y)|^2d\mu(y)d\mu(x)\big)^{1/2}.\]
	Define the Besov space $B_\sigma^{2,2}(K)$ to be $$B_{\sigma}^{2,2}(K)=\{f\in L^2(K): [f]_{B_\sigma^{2,2}(K)}<\infty\}$$ with the norm   $$\|f\|_{B_\sigma^{2,2}(K)}=\|f\|_{L^2(K)}+[f]_{B_\sigma^{2,2}(K)}.$$
\end{definition}

\begin{theorem}\label{thm48}
	For $0<\sigma<1$, we have $H^\sigma(K)=B^{2,2}_\sigma(K)$ with $\|\cdot\|_{H^\sigma(K)}\asymp \|\cdot\|_{B^{2,2}_\sigma(K)}$.
\end{theorem}

Lastly, comparing the spaces $\Lambda^{2,2}_\sigma(K)$ in terms of vertex graph approximations of $K$, we will provide some other Besov type characterizations of $H^\sigma(K)$ based on the cell graph approximations of $K$. See Theorem \ref{thm61}, Proposition \ref{prop67} and Theorem \ref{thm69}.

\section{Haar series expansion and cell graph representation}

For $m\geq 1$, $\Lambda_m$ provides a nested partition of $K$ consisting of cells with comparable diameters. Accordingly, there is a natural cell graph associated with $\Lambda_m$ inherited from the topology of $K$.

\begin{definition}\label{def51}
	For $m\geq 1$, call the graph $G_{c,m}=(\Lambda_m,E_{c,m})$ with vertex set $\Lambda_m$ and  edge set $E_{c,m}$ a level-$m$ cell graph, where $\{w,w'\}\in E_{c,m}$ if and only if $F_wK\cap F_{w'}K\neq \emptyset$.
\end{definition}

In \cite{cellgraph1}, Strichartz introduced the notions of \textit{cell graphs} and \textit{cell graph energies} for the \textit{Sierpinski gasket} $\SG$, and provided an equivalent definition of Laplacians on $\SG$, instead of using vertex graphs and vertex graph energies as we usually did in fractal analysis. See {\cite{cellgraph3,cellgraph2}} for some interesting works using related considerations.

In this section, basing on cell decomposition of functions, we will introduce a class of Besov type function spaces $\Gamma^\sigma(K)$ on $K$, which will be proved later to identify with $H^{\sigma}(K)$ when $0\leq\sigma<1$.
Recall that for $\sigma\geq 0$, we always write $\lambda=\lambda(\sigma)=r^{(d_H-\sigma d_W)/2}$ for short, where $r=\min_{i=1}^N r_i$ and $d_W=1+d_H$.

Firstly, let's look at decompositions of $L^2$ functions on $K$ by only using the nested structure provided by $\{\Lambda_m\}_{m\geq 0}$.

\subsection{Haar series expansion}
\begin{definition}\label{def52}Let $f\in L^2(K)$.
	
	(a). For any $w\in W_*$, define the average value of $f$ on $F_wK$ by \[A_w(f)=\frac{1}{\mu(F_wK)}\int_{F_wK}fd\mu.\]
	In particular, $A_\emptyset(f)=\int_{K}fd\mu$.
	
	(b). For $m\geq 1$, define the space of $m$-Haar functions
	\[\tilde{J}_m=\{\tilde{f}_m=\sum_{w\in \Lambda_m} c_w\chi_{_{F_wK}}: c_w\in \mathbb{R}, A_{w'}(\tilde{f}_m)=0,\forall w'\in \Lambda_{m-1}\}.\]
\end{definition}

It is easy to see that as subspaces of $L^2(K)$, $\tilde{J}_m\bot \tilde{J}_{m'},\forall m\neq m'$, and $$L^2(K)=constants\oplus\big(\oplus_{m=1}^\infty\tilde{J}_m\big),$$ where we use the standard inner product $(f,g)=\int_K f\cdot gd\mu$ for $f,g\in L^2(K)$. Thus for any $f\in L^2(K)$, there is a unique expansion
\[f=C+\sum_{m=1}^{\infty}\tilde{f}_m,\text{ with }\tilde{f}_m\in \tilde{J}_m.\]
The following subspaces of $L^2(K)$ are defined with weighted Haar series.

\begin{definition}\label{def53}
	For $\sigma\geq 0$, define
	\[\begin{aligned}
	\tilde{\Gamma}^\sigma(K)=\{f\in L^2(K):f=C+\sum_{m=1}^\infty\tilde{f}_m, \text{ with }\tilde{f}_m\in \tilde{J}_m, \\\text{ satisfying }\sum_{m=1}^\infty r^{-md_H}\lambda^{2m}\|\tilde{f}_m\|^2_{L^2(K)}<\infty\},
	\end{aligned}\]
	with the  norm defined by
	\[\|f\|_{\tilde{\Gamma}^\sigma(K)}=\big(|C|^2+\sum_{m=1}^\infty r^{-md_H}\lambda^{2m}\|\tilde{f}_m\|^2_{L^2(K)}\big)^{1/2}.\]
	In particular, $\tilde{\Gamma}^0(K)=L^2(K)$.
\end{definition}

\subsection{Cell graph representation}
Haar series expansion of functions  gives an easy description of functions in $L^2(K)$. However, the drawback is that the information related to the topology of $K$ is lost. To save this, we need to use the cell graphs $\{G_{c,m}\}_{m\geq 1}$.

In particular, we are interested in the cell graph energies.

\begin{definition}\label{def54} Let $m\geq 1$.
	
	(a). Define the cell graph difference operator $D_m:L^2(K)\to l(E_{c,m})$ as
	\[D_m f(\{w,w'\})=A_{w}(f)-A_{w'}(f), \quad\forall f\in L^2(K)\text{ and }\{w,w'\}\in E_{c,m}.\]
	
	(b). For $f\in L^2(K)$, define its level-$m$ cell graph energy as $\|D_mf\|^2_{l^2(E_{c,m})}$.
\end{definition}
\noindent\textbf{Remark.} In Definition \ref{def54}, we only give one value of $D_m f(\{w,w'\})$ for each edge $\{w,w'\}$, but there are two choices, say $A_{w}(f)-A_{w'}(f)$ and $A_{w'}(f)-A_{w}(f)$. However, it doesn't matter since what we only care about is the energy, and we just need a fixed choice for each edge.

Based on the cell graph energies, we define the following Besov type spaces $\Gamma^\sigma(K)$ for $\sigma\geq 0$.

\begin{definition}\label{def55}	
	For $\sigma\geq 0$, define the space
	\[\begin{aligned}
	\Gamma^\sigma(K)=\{f\in L^2(K):\sum_{m=1}^{\infty} \lambda^{2m} \|D_mf\|^2_{l^2(E_{c,m})}<\infty\},
	\end{aligned}\]
	with the norm defined by $$\|f\|_{\Gamma^\sigma(K)}=\big(\|f\|^2_{L^2(K)}+\sum_{m=1}^{\infty} \lambda^{2m}\|D_mf\|^2_{l^2(E_{c,m})}\big)^{1/2}.$$
	
\end{definition}

There is a close relation between $\Gamma^\sigma(K)$ and $\tilde{\Gamma}^\sigma(K)$. To illustrate this, we introduce following notations.

\begin{definition}\label{def56}
	(a). For each $w\in W_*$, define $W_w=\{w'\in W_*: F_{w'}K\subset F_wK\}$.
	
	(b). For $m\geq 1$, denote $\tilde{G}_{c,m}=(\Lambda_m,\tilde{E}_{c,m})$, where $\{w,w'\}\in\tilde{E}_{c,m}$ if and only if $\{w,w'\}\in E_{c,m}$ and both $w,w'$ belong to some $W_{w''}$ with $w''\in \Lambda_{m-1}$.
	
	(c). For $m\geq 1$, define the operator $\tilde{D}_m:L^2(K)\to l(\tilde{E}_{c,m})$ as
	\[\tilde{D}_m f(\{w,w'\})=A_{w}(f)-A_{w'}(f), \quad\forall f\in L^2(K)\text{ and }\{w,w'\}\in \tilde{E}_{c,m}.\]
\end{definition}

This definition is similar to Definition \ref{def51} and \ref{def54}. Clearly, $\tilde{E}_{c,m}\subset E_{c,m}$ and $\|\tilde{D}_mf\|_{l^2(\tilde{E}_{c,m})}\leq\|D_mf\|_{l^2(E_{c,m})}, \forall f\in L^2(K)$. Analogous to Definition \ref{def55}, we have the following characterization of $\tilde{\Gamma}^\sigma(K)$ for $\sigma\geq 0$.

\begin{proposition}\label{prop57}
	For $\sigma\geq 0$, we have
	\[\tilde{\Gamma}^\sigma(K)=\{f\in L^2(K):\sum_{m=1}^\infty \lambda^{2m}\|\tilde{D}_mf\|^2_{l^2(\tilde{E}_{c,m})}<\infty\},\]
	with
	\[\|f\|_{\tilde{\Gamma}^\sigma(K)}\asymp \big(\|f\|_{L^2(K)}^2+\sum_{m=1}^\infty \lambda^{2m}\|\tilde{D}_mf\|^2_{l^2(\tilde{E}_{c,m})}\big)^{1/2}.\]
\end{proposition}

\textit{Proof.} It is direct to see that
\[\tilde{D}_m\tilde{f}_{l}=0,\text{ for any } l\neq m,\tilde{f}_{l}\in \tilde{J}_l,\]
from Definition \ref{def52} and \ref{def56}. Thus, for $f=C+\sum_{l=1}^\infty \tilde{f}_l$ with $\tilde{f}_l\in \tilde{J}_l$, we have
\[\tilde{D}_mf=\tilde{D}_mC+\sum_{l=1}^\infty \tilde{D}_m\tilde{f}_l=\tilde{D}_m\tilde{f}_m.\]
The proposition follows from the observation that $\|\tilde{D}_m\tilde{f}_m\|_{l^2(\tilde{E}_{c,m})}\asymp r^{-md_H/2}\|\tilde{f}_m\|_{L^2(K)}$.\hfill$\square$
\vspace{0.2cm}

Immediately, by Proposition \ref{prop57},  for $\sigma\geq 0$, we have $\Gamma^\sigma(K)\subset\tilde{\Gamma}^\sigma(K)$, with $\|f\|_{\tilde{\Gamma}^\sigma(K)}\lesssim \|f\|_{\Gamma^\sigma(K)}$.
It is of interest to see how much information is added by introducing the cell graph structure {in the definition of $\Gamma_{\sigma}(K)$ compared} with $\tilde{\Gamma}^\sigma(K)$. We have a pair of theorems to answer this question.

\begin{theorem}\label{thm58}
	For $0\leq\sigma<\frac{d_S}{2}$, we have $\Gamma^\sigma(K)=\tilde{\Gamma}^\sigma(K)$ with $\|\cdot\|_{\Gamma^\sigma(K)}\asymp\|\cdot\|_{\tilde{\Gamma}^\sigma(K)}$. {In particular, $\Gamma^0(K)=L^2(K)$ with $\|\cdot\|_{\Gamma^0(K)}\asymp\|\cdot\|_{L^2(K)}$.}
\end{theorem}

\begin{theorem}\label{thm59}
	For $\sigma>\frac{d_S}{2}$, we have $\Gamma^\sigma(K)=\tilde{\Gamma}^\sigma(K)\cap C(K)$, and $\Gamma^\sigma(K)$ is a closed subspace of $\tilde{\Gamma}^\sigma(K)$. In particular, for $\sigma\geq 1$, we have
	\[\tilde{\Gamma}^\sigma(K)\cap C(K)=constants.\]
\end{theorem}

Theorem \ref{thm58} gives an equivalent characterization of $L^2(K)$(take $\sigma=0$), which will play a key role in the later proof of $H^\sigma(K)=\Gamma^\sigma(K)$, $0\leq \sigma<1$, while Theorem \ref{thm59} is an interesting observation from the full characterization of $H^\sigma(K)$. We prove Theorem \ref{thm58} in this section, and postpone the proof of Theorem \ref{thm59} until Section 7.

The proof of Theorem \ref{thm58} relies on the following two lemmas.

\begin{lemma}\label{lemma510}
	For  $l\geq 1$, $\tilde{f}_l\in \tilde{J}_l$, we have $\|D_l\tilde{f}_l\|_{l^2(E_{c,l})}\asymp r^{-ld_H/2}\|\tilde{f}_l\|_{L^2(K)}$.
\end{lemma}

\textit{Proof. } Obviously, we have
$r^{-ld_H/2}\|\tilde{f}_l\|_{L^2(K)}\asymp \|\tilde{D}_l\tilde{f}_l\|_{l^2(\tilde{E}_{c,l})}\leq \|D_l\tilde{f}_l\|_{l^2(E_{c,l})}.$
Conversely, the estimate $\|D_l\tilde{f}_l\|_{l^2(E_{c,l})}\lesssim r^{-ld_H/2}\|\tilde{f}_l\|_{L^2(K)}$ is also clear, since
\[\begin{aligned}
\|D_l\tilde{f}_l\|^2_{l^2(E_{c,l})}&=\sum_{\{w,w'\}\in E_{c,l}}\big(A_{w}(\tilde{f}_l)-A_{w'}(\tilde{f}_l)\big)^2\\
&\lesssim \sum_{\{w,w'\}\in E_{c,l}}\big(A^2_{w}(\tilde{f}_l)+A^2_{w'}(\tilde{f}_l)\big)\lesssim \sum_{w\in \Lambda_l}A^2_w(\tilde{f}_l)\asymp r^{-ld_H}\|\tilde{f}_l\|^2_{L^2(K)},
\end{aligned}\]
where we use the fact that the number of cells in $\Lambda_l$ neighboring $F_wK$ is bounded by $\#V_0\#\mathcal{C}$ in the second `$\lesssim$'. \hfill$\square$

\begin{lemma}\label{lemma511}
	(a). For any $\tilde{f}_l\in \tilde{J}_l$ and $l>m$, we have
	$D_m\tilde{f}_l=0.$
	
	(b). For any $\tilde{f}_l\in \tilde{J}_l$ and $l\leq m$, we have
	$\|D_m\tilde{f}_l\|_{l^2(E_{c,m})}\asymp \|D_l\tilde{f}_l\|_{l^2(E_{c,l})}$.
\end{lemma}
\textit{Proof.} (a) is obvious, so we only need to prove (b). In fact, since
\[\begin{aligned}
\|D_m\tilde{f}_l\|^2_{l^2(E_{c,m})}&=\sum_{\{w,w'\}\in E_{c,m}}\big(A_w(\tilde{f}_l)-A_{w'}(\tilde{f}_l)\big)^2\\
&=\sum_{\{v,v'\}\in E_{c,l}} \sum_{\substack{\{w,w'\}\in E_{c,m},\\w\in W_v,w'\in W_{v'}}}\big(A_w(\tilde{f}_l)-A_{w'}(\tilde{f}_l)\big)^2\\
&=\sum_{\{v,v'\}\in E_{c,l}}\sum_{\substack{\{w,w'\}\in E_{c,m},\\w\in W_v,w'\in W_{v'}}}\big(A_v(\tilde{f}_l)-A_{v'}(\tilde{f}_l)\big)^2,
\end{aligned}\]
and
$1\leq \#\big\{\{w,w'\}\in E_{c,m}:w\in W_v,w'\in W_{v'}\big\}\leq \#V_0(\#\mathcal{C})^2$, the estimate follows.
\hfill$\square$
\vspace{0.2cm}

\textit{Proof of Theorem \ref{thm58}.} By Proposition \ref{prop57}, we already have $\Gamma^\sigma(K)\subset\tilde{\Gamma}^\sigma(K)$, with $\|f\|_{\tilde{\Gamma}^\sigma(K)}\lesssim \|f\|_{\Gamma^\sigma(K)}$. It remains to show the other direction.

For $f=C+\sum_{l=1}^\infty \tilde{f}_l$ in $\tilde{\Gamma}^\sigma(K)$, using Lemma \ref{lemma511}, we get the estimate that
\[\begin{aligned}
\big(\sum_{m=1}^\infty\lambda^{2m}\|D_mf\|^2_{l^2(E_{c,m})}\big)^{1/2}&=\big\|\lambda^m\|D_mf\|_{l^2(E_{c,m})}\big\|_{l^2}=\big\|\lambda^m\|\sum_{l=1}^m D_m\tilde{f}_l\|_{l^2(E_{c,m})}\big\|_{l^2}\\
&\leq\big\|\lambda^m\sum_{l=1}^m\|D_m\tilde{f}_l\|_{l^2(E_{c,m})}\big\|_{l^2}\lesssim \big\|\lambda^m\sum_{l=1}^m\|D_l\tilde{f}_l\|_{l^2(E_{c,l})}\big\|_{l^2}\\
&=\big\|\sum_{l=0}^{m-1}\lambda^l\cdot\lambda^{m-l}\|D_{m-l}\tilde{f}_{m-l}\|_{l^2(E_{c,m-l})}\big\|_{l^2}\\
&\leq (\sum_{l=0}^\infty \lambda^l)\cdot\big\|\lambda^m\|D_m\tilde{f}_m\|_{l^2(E_{c,m})}\big\|_{l^2},
\end{aligned}\]
where we use Minkowski inequality and the fact that $\lambda<1$ in the last inequality.
Then applying Lemma \ref{lemma510}, we get $\|f\|_{\Gamma^\sigma(K)}\lesssim \|f\|_{\tilde{\Gamma}^\sigma(K)}$.\hfill$\square$

\section{Smoothed Haar functions}
In this section, we would like to establish weighted expansions of functions in $\Gamma^\sigma(K)$, analogous to Definition $\ref{def52}$ of $\tilde{\Gamma}^\sigma(K)$. Using this, we will show the Sobolev space $H^\sigma(K)$ is identical with $\Gamma^\sigma(K)$ when $0\leq\sigma<1$.

\begin{theorem}\label{thm61}
	For $0\leq\sigma<1$, $H^\sigma(K)=\Gamma^\sigma(K)$ with $\|\cdot\|_{H^\sigma(K)}\asymp \|\cdot\|_{\Gamma^\sigma(K)}$.
\end{theorem}

We will prove Theorem \ref{thm61} in the following three subsections. As an application of Theorem \ref{thm61} and the weighted expansions of functions in $\Gamma^\sigma(K)$, we will prove Theorem \ref{thm48}, i.e. $H^\sigma(K)=B^{2,2}_\sigma(K)$ for $0<\sigma<1$ in the last subsection. The following \textit{smoothed Haar functions} will play a key role.

\begin{definition}\label{def62}
	For $m\geq 1$, for any $\tilde{f}_m\in \tilde{J}_m$, define $S_m \tilde{f}_m$ to be the unique function in $dom\mathcal{E}$ such that
	\begin{equation}\label{eqn61}
	A_w(S_m\tilde{f}_m)=A_w(\tilde{f}_m), \quad\forall w\in \Lambda_m,
	\end{equation}
	and
	\begin{equation}\label{eqn62}
	\mathcal{E}(S_m\tilde{f}_m)=\min\{\mathcal{E}(f):f\in dom\mathcal{E}, A_w(f)=A_w(\tilde{f}_m),\forall w\in \Lambda_m\},
	\end{equation}
	call it a $m$-smoothed Haar function; Define the space of $m$-smoothed Haar functions by $J_m=S_m\tilde{J}_m$.
\end{definition}

\subsection{A decomposition of $dom{\mathcal{E}}$}
In this subsection, we will explore some properties of smoothed Haar functions.

First, we have an easy observation of the orthogonality of smoothed Haar function spaces as following. {We write $X\perp Y$ in $dom\mathcal{E}$ if $\mathcal{E}(f,g)=0$ for any $f\in X$ and $g\in Y$.}

\begin{lemma}\label{lemma63}
	For any $m\neq m'\geq 1$, we have $J_m\bot J_{m'}$ in dom$\mathcal{E}$.
\end{lemma}
\textit{Proof.} Without loss of generality, we assume $m<m'$. By (\ref{eqn62}), for any $f_m\in \tilde{J}_m$, and any $f\in dom\mathcal{E}$ such that $A_w(f)=0,\forall w\in \Lambda_m$, by a variational argument, we have
\[\mathcal{E}(S_m\tilde{f}_m,f)=0.\]
Then clearly
\[\mathcal{E}(S_m\tilde{f}_m,S_{m'}\tilde{f}_{m'})=0, \forall \tilde{f}_m\in \tilde{J}_m,\tilde{f}_{m'}\in \tilde{J}_{m'},\]
since $A_w(S_{m'}\tilde{f}_{m'})=0, \forall w\in \Lambda_m$ by (\ref{eqn61}).\hfill$\square$
\vspace{0.2cm}

Next, we need to estimate the energies of smoothed Haar functions. The following lemma provides a lower bound estimate.

\begin{lemma}\label{lemma64}
	For  $m\geq 1$ and $f\in dom\mathcal{E}$, we have
	\[\|D_mf\|^2_{l^2(E_{c,m})}\lesssim r^m\mathcal{E}(f).\]
\end{lemma}
\textit{Proof.} For each $\{w,w'\}\in E_{c,m}$, let $p\in F_wK\cap F_{w'}K$, by Morrey-Sobolev's inequality, we have the estimate that
\[\begin{aligned}
\big(A_w(f)-A_{w'}(f)\big)^2&\lesssim \big(f(p)-A_w(f)\big)^{2}+ \big(f(p)-A_{w'}(f)\big)^{2}\\
&\lesssim r^m\big(\mathcal{E}_{F_wK}(f)+\mathcal{E}_{F_{w'}K}(f)\big).
\end{aligned}\]
The lemma follows by summing up the estimates over all edges in $E_{c,m}$. \hfill$\square$
\vspace{0.2cm}

The upper bound estimate is due to the following lemma.
\begin{lemma}\label{lemma65}
	For $m\geq 1$ and $f\in constants\oplus(\oplus_{l=1}^m J_l)$, we have \[r^m\mathcal{E}(f)\lesssim \|D_mf\|_{l^2(E_{c,m})}^2.\]
\end{lemma}
\textit{Proof.} Noticing that
\[\mathcal{E}(f)=\inf\{\mathcal{E}(g):g\in dom\mathcal{E}, A_w(g)=A_w(f),\forall w\in \Lambda_m\},\]
we only need to construct a function $g\in dom\mathcal{E}$  with $A_w(g)=A_w(f),\forall w\in \Lambda_m$, such that $r^m\mathcal{E}(g)\lesssim \|D_mg\|_{l^2(E_{c,m})}^2=\|D_mf\|_{l^2(E_{c,m})}^2$.

This can be done as follows. First, define a piecewise harmonic function $g_0$ such that for each $x\in \bigcup_{w\in\Lambda_m} F_wV_0$,
\[g_0(x)=\frac{1}{\# \{w\in \Lambda_m: x\in F_wK\}}\sum_{w\in \Lambda_m, x\in F_wK}A_w(f),\]
and $g_0$ is harmonic in $F_wK$ for each $w\in \Lambda_m$.
Next, choose $g'\in dom_0\mathcal{E}$ such that $A_0(g')=1$. Define
\[
g=g_0+\sum_{w\in\Lambda_m}\big(A_w(f)-A_w(g_0)\big)g'\circ F_w^{-1}.
\]
For each $w\in \Lambda_m$, we have $A_w(g)=A_w(f)$ and
\[\begin{aligned}
r_w\mathcal{E}_{F_wK}(g)&=\mathcal{E}(g\circ F_w)\lesssim\mathcal{E}(g_0\circ F_w)+\big(A_w(f)-A_w(g_0)\big)^2\mathcal{E}(g')\\&\lesssim \sum_{w':\{w,w'\}\in E_{c,m}}\big(A_w(f)-A_{w'}(f)\big)^2.
\end{aligned}\]
Summing up the inequalities over $w\in\Lambda_m$, we see that $g$ satisfies the required estimate.\hfill$\square$
\vspace{0.2cm}

Combining Lemma \ref{lemma64} and  \ref{lemma65}, we get
\begin{lemma}\label{lemma66}
	For $m\geq 1$ and  $f\in constants\oplus(\oplus_{l=1}^m J_l)$, we have
	\[\|D_mf\|^2_{l^2(E_{c,m})}\asymp r^m\mathcal{E}(f).\]
\end{lemma}

The following proposition is an important corollary of the above lemma.

\begin{proposition}\label{prop67}
	dom$\mathcal{E}=constants\oplus(\oplus_{m=1}^\infty J_m)$. In addition, for any $f=C+\sum_{m=1}^\infty f_m,$ with $f_m\in J_m$, we have
	\[|C|^2+\mathcal{E}(f)\asymp \|f\|^2_{L^2(K)}+\sum_{m=1}^{\infty}r^{-m}\|D_mf_m\|^2_{l^2(E_{c,m})}.\]
\end{proposition}
\textit{Proof.} The second part of the theorem is direct by Lemma \ref{lemma63} and  \ref{lemma66}. We only need to show that $dom\mathcal{E}=constants\oplus\big(\oplus_{m=1}^\infty J_m)$.

In fact, for any $f\in dom\mathcal{E}$, let $C=A_{\emptyset}(f)$, we can inductively find a sequence $\{f_m\}_{m=1}^\infty$ such that $\forall m\geq 1$,
\[\begin{cases}
f_m\in J_m,\\
A_w(C+\sum_{l=1}^{m}f_l)=A_w(f),\forall w\in \Lambda_m.
\end{cases}\]
Then clearly we have $\mathcal{E}(\sum_{l=1}^{m}f_l)\leq \mathcal{E}(f),\forall m\geq 1$. As a result,
\[\|\sum_{l=1}^{m}f_l\|_{L^\infty(K)}^2\lesssim  \mathcal{E}(\sum_{l=1}^{m}f_l)\leq \mathcal{E}(f).\]
Thus for any $m\geq 1$ and $w\in \Lambda_m$, we have
\[A_w(C+\sum_{l=1}^{\infty}f_l)=\lim_{m'\to\infty} A_w(C+\sum_{l=1}^{m'}f_l)=A_w(C+\sum_{l=1}^{m}f_l)=A_w(f),\]
since $A_w(f_{l})=0,\forall l>m$. This shows that $f=C+\sum_{m=1}^{\infty}f_m$. \hfill$\square$

\subsection{A decomposition of $\Gamma^\sigma(K)$ with $0\leq \sigma<1$}
The benefit of the smoothed Haar functions comes from two aspects. First, they keep features from Haar functions. Second, they are ``smooth''. Below we give a lemma which follows from what we have discussed in the last subsection.

\begin{lemma}\label{lemma68}
	(a). For any $m'\geq m\geq 1$ and any $f\in constants\oplus(\oplus_{l=1}^m J_l)$, we have
	\[\|D_{m'}f\|_{l^2(E_{c,m'})}\lesssim r^{\frac{m'-m}{2}}\|D_mf\|_{l^2(E_{c,m})}.\]
	
	(b). For any $m'>m\geq 1$ and any $f\in J_{m'}$, we have $D_{m}f=0$.
\end{lemma}

\textit{Proof.} (a) is a simple consequence of Lemma \ref{lemma64} and  \ref{lemma66}. (b) comes naturally from Definition \ref{def62}.\hfill$\square$
\vspace{0.2cm}

\begin{theorem}\label{thm69}
	Let $0\leq\sigma<1$ and  $\lambda:=\lambda(\sigma)=r^{(d_H-\sigma d_W)/2}$. We have
	
	(a). For any $f\in \Gamma^\sigma(K)$, there is a unique sequence of functions $\{f_m\}_{m=1}^\infty$ with $f_m\in J_m$ and a constant $C$ such that $f=C+\sum_{m=1}^\infty f_m$.
	
	(b). For any series $f=C+\sum_{m=1}^{\infty}f_m$, $f\in \Gamma^\sigma(K)$ if and only if
	\[\sum_{m=1}^{\infty} \lambda^{2m}\|D_mf_m\|^2_{l^2(E_{c,m})}<\infty.\]
	In addition,
	\[\|f\|_{\Gamma^\sigma(K)}\asymp \big(|C|^2+\sum_{m=1}^{\infty}\lambda^{2m}\|D_mf_m\|^2_{l^2(E_{c,m})}\big)^{1/2}.\]
\end{theorem}

\textit{Proof.} (a) and (b) are linked with each other, so we will prove them at the same time.

First, let $\{f_l\}_{l=1}^\infty$ be a sequence such that $f_l\in J_l,\forall l\geq 1$, and
\[\big\|\lambda^l\|D_lf_l\|_{l^2(E_{c,l})}\big\|^2_{l^2}=\sum_{l=1}^{\infty}\lambda^{2l}\|D_lf_l\|^2_{l^2(E_{c,l})}<\infty.\]
Let $C$ be a constant. We first need to show that the series $C+\sum_{l=1}^\infty f_l$ converges in $L^2(K)$, which is enough by showing that $\{C+\sum_{l=1}^nf_l\}_{n=0}^\infty$ is a Cauchy sequence. In fact, by using Theorem \ref{thm58}, for any $n'\geq n$, we see that
$$\big\|\sum_{l=n}^{n'} f_l\big\|_{L^2(K)}\lesssim  \Big(\sum_{m=1}^\infty \big\|r^{md_H/2}D_m\big(\sum_{l=n}^{n'} f_l\big)\big\|_{l^2(E_{c,m})}^2\Big)^{1/2},$$
and thus
\begin{equation}\label{eqn63}
\begin{aligned}
\big\|\sum_{l=n}^{n'} f_l\big\|_{L^2(K)}&\lesssim \big\|\lambda^m\sum_{l=n}^{n'} \|D_mf_l\|_{l^2(E_{c,m})}\big\|_{l^2}\\
&\lesssim \big\|\lambda^m\sum_{l=n}^{m} r^{(m-l)/2}\|D_lf_l\|_{l^2(E_{c,l})}\big\|_{l^2}\\
&=\big\|\lambda^m\sum_{l=0}^{m-n} r^{l/2}\|D_{m-l}f_{m-l}\|_{l^2(E_{c,m-l})}\big\|_{l^2}\\
&= \big\|\sum_{l=0}^{m-n}\big(r^{1/2}\lambda\big)^l\cdot \lambda^{m-l}\|D_{m-l}f_{m-l}\|_{l^2(E_{c,m-l})}\big\|_{l^2}\\
&\leq \big\|\chi_{m\geq n}\lambda^m\|D_mf_m\|_{l^2(E_{c,m})}\big\|_{l^2}\sum_{l=0}^{\infty}\big(r^{1/2}\lambda\big)^l\\&\lesssim \big\|\chi_{m\geq n}\lambda^m\|D_mf_m\|_{l^2(E_{c,m})}\big\|_{l^2},
\end{aligned}
\end{equation}
where we used Lemma \ref{lemma68} in the second inequality, Minkowski inequality in the last but one inequality, and the fact that $r^{1/2}\lambda<1$. This gives that \[\lim_{n,n'\to\infty}\big\|\sum_{l=n}^{n'} f_l\big\|_{L^2(K)}=0,\]
and so that $C+\sum_{l=0}^\infty f_l$ is well defined in $L^2(K)$. Let $f=C+\sum_{l=0}^\infty f_l$. As a consequence, for any $m\geq1$, $w\in \Lambda_m$, we get
\begin{equation}\label{eqn64}
A_w(f)=\lim_{m'\to\infty}A_w(C+\sum_{l=1}^{m'}f_l)=A_w(C+\sum_{l=1}^{m}f_l),
\end{equation}
since $A_w$ is a bounded functional on $L^2(K)$. Thus,
\[
D_m(f)=\sum_{l=1}^m D_mf_l.
\]
Using the above equality and by a similar process as estimate (\ref{eqn63}), we get the estimate that
\begin{equation}\label{eqn65}
\|f\|_{\Gamma^\sigma(K)}\lesssim \big(|C|^2+\sum_{m=1}^\infty{ \lambda^{2m}}\|D_mf_m\|^2_{l^2(E_{c,m})}\big)^{1/2}.
\end{equation}

Next, let $f\in \Gamma^\sigma(K)$. Let $C=A_\emptyset (f)$, we can inductively find a sequence $\{f_m\}_{m=1}^\infty$ such that
\[\begin{cases}
f_m\in J_m,&\forall m\geq 1,\\
A_w(C+\sum_{l=1}^mf_l)=A_w(f),&\forall w\in \Lambda_m, m\geq 0.
\end{cases}\]
For simplicity, we write $g_m=C+\sum_{l=1}^{m}f_l\in constants\oplus (\oplus_{l=1}^mJ_l)$. Then by Lemma \ref{lemma68} (a), for $m\geq 1$, we have
\[\|D_{m+1}g_{m}\|_{l^2(E_{c,m+1})}\lesssim \|D_{m}g_{m}\|_{l^2(E_{c,m})}=\|D_{m}f\|_{l^2(E_{c,m})}.\]
As a result,
\[\begin{aligned}
\|D_{m+1}f_{m+1}\|_{l^2(E_{c,m+1})}&\leq \|D_{m+1}g_{m+1}\|_{l^2(E_{c,m+1})}+\|D_{m+1}g_{m}\|_{l^2(E_{c,m+1})}\\
&\lesssim \|D_{m+1}f\|_{l^2(E_{c,m+1})}+\|D_{m}f\|_{l^2(E_{c,m})}.
\end{aligned}\]
In particular, $\|D_1f_1\|_{l^2(E_{c,1})}=\|D_1f\|_{l^2(E_{c,1})}$.
Thus we get
\begin{equation}\label{eqn66}
\sum_{m=1}^{\infty} \lambda^{2m}\|D_mf_m\|^2_{l^2(E_{c,m})}\lesssim \sum_{m=1}^{\infty} \lambda^{2m}\|D_mf\|^2_{l^2(E_{c,m})}.
\end{equation}
By the first part of the discussion, we can see that $f=C+\sum_{m=1}^{\infty}f_m$, by using (\ref{eqn64}).

Thus we have proved the theorem. The estimates in (b) come from (\ref{eqn65}) and (\ref{eqn66}).\hfill$\square$

\subsection{Proof of  $H^\sigma(K)=\Gamma^\sigma(K)$ with $0\leq \sigma<1$}
We have decomposed both $dom\mathcal{E}$ and $\Gamma^\sigma(K)$ with $0\leq \sigma<1$ into summations of smoothed Haar functions in Proposition \ref{prop67} and Theorem \ref{thm69} respectively.

Recall the fact that $\|D_mS_m\tilde{f}_m\|_{l^2(E_{c,m})}=\|D_m\tilde{f}_m\|_{l^2(E_{c,m})}\asymp r^{-md_H/2}\|\tilde{f}_m\|_{L^2(K)}$ by Lemma \ref{lemma510}, and recall Definition \ref{def53} of $\tilde{\Gamma}^\sigma(K)$. The following lemma holds as an immediate corollary of Proposition \ref{prop67} and Theorem \ref{thm69}.

\begin{lemma}\label{lemma610}
	Given any function $f=C+\sum_{m=1}^\infty \tilde{f}_m\in L^2(K)$ with $\tilde{f}_m\in \tilde{J}_m$, let $Sf$ be
	\[Sf=C+\sum_{m=1}^\infty S_m\tilde{f}_m.\]
	Then $S$ is a homeomorphism from $\tilde{\Gamma}^\sigma(K)$ onto $\Gamma^\sigma(K)$, $\forall 0\leq \sigma<1$, and from $\tilde{\Gamma}^1(K)$ onto dom$\mathcal{E}$.
\end{lemma}

Since both $\tilde{\Gamma}^\sigma(K)$ and $H^\sigma(K)$ are stable under complex interpolation, the following lemma follows from Lemma \ref{lemma610} easily.

\begin{lemma}\label{lemma611}
	For $0\leq\sigma\leq 1$, $S$ is a linear homeomorphism from $\tilde{\Gamma}^\sigma(K)$ onto $H^\sigma(K)$.
\end{lemma}
\textit{Proof.} Notice that $H^0(K)=\Gamma^0(K)$ by Theorem \ref{thm58},  and  $H^1(K)=dom\mathcal{E}$ by Corollary \ref{coro33}. So the lemma is true for $\sigma=0$ and $\sigma=1$ by Lemma \ref{lemma610}, and thus for $0<\sigma<1$, by complex interpolation, since $\tilde{\Gamma}^\sigma(K)=[\tilde{\Gamma}^0(K),\tilde{\Gamma}^1(K)]_{\sigma}$ and $H^\sigma(K)=[H^0(K),H^1(K)]_{\sigma}$.\hfill$\square$
\vspace{0.2cm}

Theorem \ref{thm61} is an immediate consequence of Lemma \ref{lemma610} and \ref{lemma611}.

\subsection{Proof of $H^\sigma(K)=B^{2,2}_{\sigma}(K)$ with $0<\sigma<1$}
As an application of Theorem \ref{thm61} and \ref{thm69}, we provide a pure analytic proof of Theorem \ref{thm48}. First, we list some easy estimates. Recall  Definition \ref{def47} for the definition of $B^{2,2}_\sigma(K)$, and still write $\lambda=\lambda(\sigma)=r^{(d_H-\sigma d_W)/2}$ for short.

\begin{lemma}\label{lemma612}
	For  $f\in L^2(K)$, $m\geq 0$, write  	 \[I_m(f)=r^{-md_H}\big(\int_K\int_{B(x,r^m)}|f(x)-f(y)|^2 d\mu(y)d\mu(x)\big)^{1/2}.\]
	Then for $\sigma>0$, we have  $[f]_{B^{2,2}_\sigma(K)}\asymp \|\lambda^mI_m(f)\|_{l^2}$.
\end{lemma}
\textit{Proof.}  It is not hard to see that
\[[f]_{B^{2,2}_\sigma(K)}\asymp \big(\sum_{m=1}^{\infty} r^{-md_H-m\sigma d_W} \int_K\int_{B(x,r^m)}|f(x)-f(y)|^2d\mu(y)d\mu(x)\big)^{1/2}.\]
The lemma follows immediately since $(r^{-md_H}\lambda^{m})^2=r^{-md_H-m\sigma d_W}$. \hfill$\square$

The following are some easy estimates.
\begin{lemma}\label{lemma613}
	Let $f_l\in J_l$ for some $l\geq 1$. For $w\in \Lambda_l$, write $M_w(f_l)=\max_{x\in F_wK} |f_l(x)|$. Then we have
	\[\big(\sum_{w\in \Lambda_l}M_w^2(f_l)\big)^{1/2}\lesssim\|D_lf_l\|_{l^2(E_{c,l})}.\]
\end{lemma}

\textit{Proof.} In fact, the lemma follows from energy estimates. For $w\in \Lambda_l$, let $x\in F_wK$ such that $|f_l(x)|=M_w(f_l)$, and let $y\in F_wK$ such that $f_l(y)=A_w(f_l)$. Then
\[\mathcal{E}_{F_w(K)}(f_l)\geq R^{-1}(x,y)\big(f_l(x)-f_l(y)\big)^2\gtrsim r^{-l} \big(f_l(x)-f_l(y)\big)^2.\]
Thus by Lemma \ref{lemma66}, we have
\[r^{-l}\sum_{w\in \Lambda_l}\big({|A_w(f_l)|}-M_w(f_l)\big)^2\lesssim \mathcal{E}(f_l)\asymp r^{-l}\|D_lf_l\|^2_{l^2(E_{c,l})}.\]
The lemma follows.\hfill$\square$
\vspace{0.2cm}

As a consequence of Lemma \ref{lemma613}, we have
\begin{lemma}\label{lemma614}
	For $l\geq 1$ and  $f_l\in J_l$, we have $\|f_l\|_{L^2(K)}\asymp r^{ld_H/2}\|D_lf_l\|_{l^2(E_{c,l})}$.
\end{lemma}
\textit{Proof.} Clearly, we have
\[r^{ld_H}\|D_lf_l\|^2_{l^2(E_{c,l})}\asymp r^{ld_H}\sum_{w\in \Lambda_l}|A_w(f_l)|^2\lesssim \|f_l\|^2_{L^2(K)}\lesssim r^{ld_H}\sum_{w\in \Lambda_l}|M_w(f_l)|^2.\]
The lemma follows immediately from Lemma \ref{lemma613}.\hfill$\square$
\vspace{0.2cm}

We have the following estimate of $I_m(f_l)$ with $f_l\in J_l$ and $m<l$.
\begin{lemma}\label{lemma615}
	$I_m(f_l)\lesssim r^{\frac{l-m}{2}d_H}\|D_lf_l\|_{l^2(E_{c,l})}$ for $f_l\in J_l$ and $0\leq m<l$.
\end{lemma}
\textit{Proof.} Note that
\[\begin{aligned}
I_m^2(f_l)&=r^{-2md_H}\int_{K}\int_{B(x,r^m)}|f_l(x)-f_l(y)|^2d\mu(y)d\mu(x)\\
&\leq 2\cdot r^{-2md_H}\int_{K}\int_{B(x,r^m)}(|f_l(x)|^2+|f_l(y)|^2)d\mu(y)d\mu(x)\\
&\lesssim r^{-md_H}\|f_l\|^2_{L^2(K)}\lesssim r^{(l-m)d_H}\|D_lf_l\|^2_{l^2(E_{c,l})},
\end{aligned}\]
where we use Lemma \ref{lemma614} in the last inequality.\hfill$\square$
\vspace{0.2cm}

The following lemma is also useful.
\begin{lemma}\label{lemma616}
	For $m\geq 0$ and $f\in dom\mathcal{E}$, we have \[I_m(f)\lesssim r^{m/2}\mathcal{E}^{1/2}(f).\]
\end{lemma}
\textit{Proof.} We only need to notice that for any $w\in \Lambda_m$, we have
\[r^{-2md_H}\int_{F_wK}\int_{B(x,r^m)} |f(x)-f(y)|^2d\mu(y)d\mu(x)\lesssim r^m\mathcal{E}_A(f),\]
where $A=\bigcup\{F_{w'}K:\exists x\in F_wK \text{ so that }F_{w'}K\cap B(x,r^m)\neq \emptyset\}$, as $|f(x)-f(y)|^2\lesssim r^m\mathcal{E}_A(f)$. The lemma follows by summing the above estimate over all $w\in\Lambda_m$, using Lemma \ref{lemma25}.\hfill$\square$
\vspace{0.2cm}

\noindent\textbf{Remark.} In\cite{GL}, Gu and Lau studied another class of Besov spaces, $B^{2,\infty}_\sigma(K)$, which includes $dom\mathcal{E}$ as a critical case. They proved that for $f\in dom\mathcal{E}$, $\sup_{m\geq 0}r^{-m}I_m^2(f)\asymp \mathcal{E}(f)$.
\vspace{0.2cm}

As a consequence of Lemma \ref{lemma66} and \ref{lemma616}, we get
\begin{lemma}\label{lemma617}
	For $1\leq m\leq m'$, and $f\in constants\oplus_{l=1}^{m'} J_l$, we have \[I_m(f)\lesssim r^{\frac{m-m'}{2}}\|D_{m'}f\|_{l^2(E_{c,m'})}.\]
\end{lemma}

\textit{Proof of Theorem \ref{thm48}}. By Theorem \ref{thm61}, it suffices to prove $\Gamma^\sigma(K)=B_{\sigma}^{2,2}(K)$.

First, let's prove  that $B_{\sigma}^{2,2}(K)\subset \Gamma^\sigma(K)$. Without loss of generality, assume that $diam(K)=\max\{R(x,y): x,y\in K\}=1$. Then for $f\in B_\sigma^{2,2}(K)$, for $m\geq 1$ and  $\{w,w'\}\in E_{c,m}$, we have
\[\begin{aligned}
|A_w(f)-A_{w'}(f)|^2&=\big|\frac{1}{\mu(F_wK)}\int_{F_wK}f(x)d\mu(x)-A_{w'}(f)\big|^2\\
&\lesssim \big(\int_{F_wK}r^{-md_H}|f(x)-A_{w'}(f)|d\mu(x)\big)^2\\
&\leq r^{-md_H}\int_{F_wK} |f(x)-A_{w'}(f)|^2d\mu(x)\\
&\leq r^{-2md_H}\int_{F_wK}\int_{F_{w'}K} |f(x)-f(y)|^2d\mu(y)d\mu(x)\\
&\leq r^{-2md_H}\int_{F_wK}\int_{B(x,2r^m )} |f(x)-f(y)|^2d\mu(y)d\mu(x),
\end{aligned}\]
where we use  Jensen's inequality in the third inequality. Summing the above estimate over all edges in $E_{c,m}$, we then get that
\[\|D_mf\|_{l^2(E_{c,m})}\lesssim I_{m+[\log 2/\log r]}(f).\]
This gives that  $f\in \Gamma^\sigma(K)$ and $\|f\|_{\Gamma^\sigma(K)}\lesssim\|f\|_{B_\sigma^{2,2}(K)}$ by Lemma \ref{lemma612}.

Next, for the other direction, recall that by Theorem \ref{thm69}, for any $f\in \Gamma^\sigma(K)$, there is a unique expansion that $f=C+\sum_{l=1}^{\infty}f_l$ with $f_l\in J_l$. For $m\geq 1$, we write $g_m=C+\sum_{l=1}^{m}f_l$. Then we have the estimate that
\[\begin{aligned}
\|\lambda^mI_m(f)\|_{l^2}&\lesssim \|\lambda^mI_m(g_m)+\lambda^m\sum_{l=m+1}^\infty I_m(f_l)\|_{l^2}\\
&\leq\|\lambda^mI_m(g_m)\|_{l^2}+\|\lambda^m\sum_{l=m+1}^\infty I_m(f_l)\|_{l^2}\\
&\lesssim \big\|\lambda^m\|D_{m}g_m\|_{l(E_{c,m})}\big\|_{l^2}+\big\|\lambda^m\sum_{l=m+1}^{\infty}r^{\frac{l-m}2d_H}\|D_lf_l\|_{l^2(E_{c,l})}\big\|_{l^2}\\
&=
\big\|\lambda^m\|D_mf\|_{l(E_{c,m})}\big\|_{l^2}+\big\|\lambda^m\sum_{l=1}^{\infty}r^{ld_H/2}\|D_{m+l}f_{m+l}\|_{l^2(E_{c,m+l})}\big\|_{l^2}\\
&\leq \big\|\lambda^m\|D_mf\|_{l(E_{c,m})}\big\|_{l^2}+\sum_{l=1}^{\infty}(r^{-d_H/2}\lambda)^{-l}\cdot\big\|\lambda^{m}\|D_{m}f_{m}\|_{l^2(E_{c,m})}\big\|_{l^2}\\
&\lesssim \|f\|_{\Gamma^\sigma(K)},
\end{aligned}\]
where we use Lemma \ref{lemma615} and \ref{lemma617} in the third inequality,  use Minkowski inequality in the fourth inequality, and use Theorem \ref{thm69} in the last inequality, noticing that $r^{-d_H/2}\lambda>1$.  This finishes the proof.\hfill$\square$
\vspace{0.2cm}

\noindent\textbf{Remark.} In fact, in the above proof, we have proved that $B_{\sigma}^{2,2}(K)\subset\Gamma^\sigma(K)$ for all $\sigma>0$. Combining this with Theorem \ref{thm59}, we immediately get that $B^{2,2}_\sigma(K)=constants$ whenever $\sigma\geq \frac 12$.

\section{Atomic decompositions: lower orders}
In Section 5 and 6, we have developed all the tools for the proof of the atomic decomposition theorem Theorem \ref{decompositionthm} for $0\leq \sigma<1$. Now we will come to the proof of this case in this section. In the next section, we then sketch the proof for $1\leq\sigma<2$. The general $\sigma\geq 0$ case then immediately follows with formula (\ref{eqn34}) by combining the results for $0\leq\sigma<1$ and $1\leq\sigma<2$ together.

As applications of the case $0\leq\sigma<1$, additionally in this section, we will prove Theorem \ref{thm46} (a), i.e.  for $\frac{d_S}{2}<\sigma<1$, $H^\sigma(K)$ is identical with another class of Besov type spaces  $\Lambda_{\sigma}^{2,2}(K)$. Also, we will prove Theorem \ref{thm59}.

As in previous sections, we still abbreviate that $\lambda=\lambda(\sigma)=r^{(d_H-\sigma d_W)/2}$ for $\sigma\geq 0$.

\subsection{Atomic decomposition: $0\leq \sigma<\frac{d_S}{2}$}

When $\sigma<\frac{d_S}{2}$, we already have $H^\sigma(K)=\tilde{\Gamma}^\sigma(K)$ by Theorem \ref{thm58} and \ref{thm61}. This immediately implies Theorem \ref{decompositionthm} for the case $0\leq \sigma<\frac{d_S}{2}$. \vspace{0.2cm}

\textit{Proof of Theorem \ref{decompositionthm} when $0\leq \sigma<\frac{d_S}{2}$.}
For each $m\geq 1$, write
\begin{equation}\label{eqn71}
\Lambda^{-}_m=\{w\in W_*: \exists w'\in\Lambda_{m-1},w''\in\Lambda_m, \text{ such that } F_{w''}K\subsetneq F_wK\subset F_{w'}K\},
\end{equation}
and we can easily check that {$W_*=\bigcup_{m=1}^\infty \Lambda^{-}_m$}.

Define
\[\tilde{f}_m=\sum_{w\in \Lambda^{-}_m}\sum_{i=1}^N a_{wi}\chi_{{F_{wi}K}},\]
Then, it is direct to check that $\tilde{f}_m\in \tilde{J}_m,\forall m\geq 1$. In addition,
\[
r^{-md_H}\lambda^{2m}\|\tilde{f}_m\|^2_{L^2(K)}=r^{-m\sigma d_W}\sum_{w\in \Lambda^-_{m}}\sum_{i=1}^N r_{wi}^{d_H}|a_{wi}|^2\asymp \sum_{w\in \Lambda^-_{m}}r_{w}^{d_H-\sigma d_W}\sum_{i=1}^N |a_{wi}|^2.
\]
The result follows from the above estimate immediately, noticing that $\|f\|_{H^\sigma(K)}\asymp \big(|C|^2+\sum_{m=1}^\infty r^{-md_H}\lambda^{2m}\|\tilde{f}_m\|^2_{L^2(K)}\big)^{1/2}$. \hfill$\square$

\subsection{Atomic decomposition: $\frac{d_S}{2}<\sigma<1$} In this part, we will prove Theorem \ref{decompositionthm} for $\frac{d_S}{2}<\sigma<1$.  Recall the definition of tent functions $\psi_x$ in Section 4.1.
\begin{definition}
For $m\geq 1$, denote $T_m$ the linear subspace spanned by $\{\psi_{F_wx}:w\in \Lambda^-_m,x\in V_1\setminus V_0\}$, where $\Lambda^-_m$ is defined in (\ref{eqn71}).
\end{definition}

When $\sigma>\frac{d_S}{2}$, it is well-known that $H^\sigma(K)\subset C(K)$. In fact, this is an immediate consequence of Theorem \ref{thm69}. Obviously, each function $f\in C(K)$ can be written uniquely as a series of tent functions
\begin{equation}\label{eqn72}
f=\varphi_0+\sum_{x\in V_*\setminus V_0} c_x\psi_x,
\end{equation}
with $\varphi_0\in \mathcal{H}_0$ and $c_x\in \mathbb{R}$. By the definition of $T_m$, immediately we see that $f$ admits a unique expansion of the form
\begin{equation}\label{eqn73}
f=\sum_{m=0}^\infty \varphi_m,
\end{equation}
with $\varphi_0\in \mathcal{H}_0$ and $\varphi_m\in T_m,\forall m\geq 1$.

The following lemma follows by a similar argument as the proof in Section 7.1.

\begin{lemma}\label{lemma72}
	Let $f\in C(K)$, and take the expansions in (\ref{eqn72}) and (\ref{eqn73}). Then, for $\sigma>\frac{d_S}{2}$, we have $\sum_{w\in W_*}r_w^{d_H-\sigma d_W}\sum_{x\in V_1\setminus V_0}|c_{F_wx}|^2<\infty$ if and only if $\sum_{m=0}^\infty r^{-m\sigma d_W}\|\varphi_m\|^2_{L^2(K)}<\infty$. In addition,
	\[\|\varphi_0\|^2_{L^2(K)}+\sum_{w\in W_*}r_w^{d_H-\sigma d_W}\sum_{x\in V_1\setminus V_0}|c_{F_wx}|^2\asymp \sum_{m=0}^\infty r^{-m\sigma d_W}\|\varphi_m\|^2_{L^2(K)}.\]
\end{lemma}

One can compare Lemma \ref{lemma72} with Theorem \ref{thm69}. In fact, we will see that there is mutual control of norms of functions in $J_m$ and $T_m$. The same idea will be used in the proof for higher orders in Section 8.

Recall $\nabla_m$ defined in Definition \ref{def42}. It is direct to check that for $m\geq 1$ and $\varphi_m\in T_m$, we always have
\begin{equation}\label{eqn74}
r^{-md_H/2}\|\varphi_m\|_{L^2(K)}\asymp\|\nabla_m\varphi_m\|_{l^2(E_{v,m})}\asymp r^{m/2}\mathcal{E}^{1/2}(\varphi_m).
\end{equation}

The following lemmas provide a mutual control of functions in $J_m$ and $T_m$.
\begin{lemma}\label{lemma73}
	Let $l\geq 1$, $f_l=\sum_{m=0}^\infty\varphi_{m,l}\in J_l$ with $\varphi_{0,l}\in \mathcal{H}_0$ and $\varphi_{m,l}\in T_m$.
	
	(a). If $m\leq l$, then $\|\nabla_m\varphi_{m,l}\|_{l^2(E_{v,m})}\lesssim \|D_lf_l\|_{l^2(E_{c,l})}$.
	
	(b). If $m>l$, then $\|\nabla_m\varphi_{m,l}\|_{l^2(E_{v,m})}\lesssim r^{(m-l)/2}\|D_lf_l\|_{l^2(E_{c,l})}$.
\end{lemma}
\textit{Proof.} Since $\varphi_{m,l}=\sum_{m'=0}^{m}\varphi_{m',l}-\sum_{m'=0}^{m-1}\varphi_{m',l}$, we have
$$\|\nabla_{m}\varphi_{m,l}\|_{l^2(E_{v,m})}\leq \big\|\nabla_m(\sum_{m'=0}^{m}\varphi_{m',l})\big\|_{l^2(E_{v,m})}+\big\|\nabla_m(\sum_{m'=0}^{m-1}\varphi_{m',l})\big\|_{l^2(E_{v,m})},$$
and thus
$$\|\nabla_{m}\varphi_{m,l}\|_{l^2(E_{v,m})}\lesssim \|\nabla_mf_l\|_{l^2(E_{v,m})}+\|\nabla_{m-1}f_l\|_{l^2(E_{v,m-1})}.$$

Then (a) follows by Lemma \ref{lemma613}, and (b) follows by Lemma \ref{lemma66}. In fact, for $m>l$, we then have
\[\quad\quad\|\nabla_m\varphi_{m,l}\|_{l^2(E_{v,m})}\lesssim r^{m/2}\mathcal{E}^{1/2}(f_l)\lesssim r^{(m-l)/2}\|D_lf_l\|_{l^2(E_{c,l})}.\quad\quad\square\]

\begin{lemma}\label{lemma74}
	Let $l\geq 0$, $\varphi_l=C_l+\sum_{m=1}^\infty f_{m,l}\in T_l$(or $\mathcal{H}_0$) with $C_l\in\mathbb{R}$ and $f_{m,l}\in J_m$.
	
	(a). If $m\leq l$, then $\|D_mf_{m,l}\|_{l^2(E_{c,m})}\lesssim r^{(l-m)d_H/2}\|\nabla_l\varphi_l\|_{l^2(E_{v,l})}$.
	
	(b). If $m>l$, then $\|D_mf_{m,l}\|_{l^2(E_{c,m})}\lesssim r^{(m-l)/2}\|\nabla_l\varphi_l\|_{l^2(E_{v,l})}$.
\end{lemma}
\textit{Proof.} (a). By Theorem \ref{thm58}(taking $\sigma=0$), and formula (\ref{eqn74}), we have $$\|D_mf_{m,l}\|_{l^2(E_{c,m})}\lesssim r^{-md_H/2}\|\varphi_l\|_{L^2(K)}\lesssim r^{(l-m)d_H/2}\|\nabla_l\varphi_l\|_{l^2(E_{v,l})}.$$

(b). By Proposition \ref{prop67} and formula (\ref{eqn74}), we have
$$\|D_mf_{m,l}\|_{l^2(E_{c,m})}\lesssim r^{m/2}\mathcal{E}^{1/2}(f_{m,l})\leq r^{m/2}\mathcal{E}^{1/2}(\varphi_l)\lesssim r^{(m-l)/2}\|\nabla_l\varphi_l\|_{l^2(E_{v,l})}.\square$$

\textit{Proof of Theorem \ref{decompositionthm} when $\frac{d_S}{2}<\sigma<1$.}  Let $f\in C(K)$.  Write $f$ with the expansions that $f=\sum_{l=0}^\infty\varphi_l$ and $f=C+\sum_{l=1}^\infty f_l$ with $\varphi_0\in \mathcal{H}_0$, $C\in\mathbb{R}$, $\varphi_l\in T_l$ and $f_l\in J_l$ for $l\geq 1$. By Lemma \ref{lemma72}, Theorem \ref{thm61}, Theorem \ref{thm69}, and using formula (\ref{eqn74}), it suffices to prove that $\sum_{l=1}^\infty\lambda^{2l}\|\nabla_l\varphi_l\|^2_{l^2(E_{v,l})}<\infty$ if and only if  $\sum_{l=1}^\infty \lambda^{2l}\|D_lf_l\|^2_{l^2(E_{c,l})}<\infty$, and
\[\|\varphi_0\|_{L^2(K)}+\sum_{l=1}^\infty\lambda^{2l}\|\nabla_l\varphi_l\|^2_{l^2(E_{v,l})}\asymp |C|^2+\sum_{l=1}^\infty \lambda^{2l}\|D_lf_l\|^2_{l^2(E_{c,l})}.\]

Assume that $\sum_{l=1}^\infty \lambda^{2l}\|\nabla_l \varphi_l\|_{l^2(E_{c,l})}^2<\infty$. For each $l\geq 0$, we expand $\varphi_l=C_l+\sum_{m=1}^\infty f_{m,l}$ with $C_l\in\mathbb{R}$ and $f_{m,l}\in J_m$. Then it is direct to check that $f_m=\sum_{l=0}^\infty f_{m,l}, \forall m\geq 1$, and thus
\[\begin{aligned}
&\big\|\lambda^m\|D_mf_m\|_{l^2(E_{c,m})}\big\|_{l^2}\leq\big\|\lambda^m\sum_{l=0}^\infty\|D_mf_{m,l}\|_{l^2(E_{c,m})}\big\|_{l^2}\\
\leq&\big\|\lambda^m\sum_{l=0}^{m-1}\|D_mf_{m,l}\|_{l^2(E_{c,m})}\big\|_{l^2}+\big\|\lambda^m\sum_{l=m}^\infty\|D_mf_{m,l}\|_{l^2(E_{c,m})}\big\|_{l^2}.
\end{aligned}\]

By Lemma \ref{lemma74} (b), using Minkowski inequality, noticing that $\lambda r^{1/2}<1$, we have
\[\begin{aligned}
\big\|\lambda^m\sum_{l=0}^{m-1}\|D_mf_{m,l}\|_{l^2(E_{c,m})}\big\|_{l^2}&\lesssim \big\|\lambda^{m}\sum_{l=0}^{m-1} r^{(m-l)/2}\|\nabla_l\varphi_l\|_{l^2(E_{v,l})}\big\|_{l^2}\\
&=\big\|\sum_{l=1}^{m} r^{l/2}\lambda^l\lambda^{m-l}\|\nabla_{m-l}\varphi_{m-l}\|_{l^2(E_{v,m-l})}\big\|_{l^2}\\
&\lesssim \big\|\lambda^m\|\nabla_m\varphi_m\|_{l^2(E_{v,m})}\big\|_{l^2}.
\end{aligned}
\]

On the other hand, by Lemma \ref{lemma74} (a), using a similar argument, we have
\[\begin{aligned}
\big\|\lambda^m\sum_{l=m}^\infty\|D_mf_{m,l}\|_{l^2(E_{c,m})}\big\|_{l^2}&\lesssim \big\|\lambda^m\sum_{l=m}^\infty r^{(l-m)d_H/2}\|\nabla_l\varphi_l\|_{l^2(E_{v,l})}\big\|_{l^2}\\
&=\|\sum_{l=0}^\infty \lambda^{-l}r^{ld_H/2}\lambda^{l+m}\|\nabla_{l+m}\varphi_{l+m}\|_{l^2(E_{v,l+m})}\big\|_{l^2}\\
&\lesssim \big\|\lambda^m\|\nabla_m\varphi_m\|_{l^2(E_{v,m})}\big\|_{l^2}.
\end{aligned}\]

Combining the above estimates, together with the observation that $|C|\leq\sum_{m=0}^\infty \|\varphi_m\|_{L^\infty(K)}$ and $\|\nabla_0\varphi_0\|_{l^2(E_{v,0})}\lesssim \|\varphi_0\|_{L^2(K)}$, we get that $$|C|^2+\sum_{m=1}^\infty \lambda^{2m}\|D_mf_m\|^2_{l^2(E_{c,m})}\lesssim \|\varphi_0\|^2_{L^2(K)}+\sum_{m=1}^\infty \lambda^{2m}\|\nabla_m\varphi_m\|^2_{l^2(E_{v,m})}.$$

The other direction follows by a similar argument by using Lemma \ref{lemma73}.\hfill$\square$

\subsection{Proof of  $H^\sigma(K)=\Lambda^{2,2}_\sigma(K)$ with $\frac{d_S}{2}<\sigma<1$} Now we come to the proof of Theorem \ref{thm46} (a). Recall Definition \ref{def43} for the Besov type spaces $\Lambda^{2,2}_\sigma(K)$.

Notice that in the last subsection we have shown that $f=\sum_{m=0}^\infty\varphi_m\in H^\sigma(K)$ if and only if $\sum_{m=1}^\infty\lambda^{2m}\|\nabla_m\varphi_m\|^2_{l^2(E_{v,m})}<\infty$, and $\|f\|_{H^\sigma(K)}\asymp \big(\|\varphi_0\|_{L^2(K)}+\sum_{m=1}^\infty\lambda^{2m}\|\nabla_m\varphi_m\|^2_{l^2(E_{v,m})}\big)^{1/2}$. Thus to prove Theorem \ref{thm46} (a), it suffices to prove a same result with $\Lambda^{2,2}_\sigma(K)$ instead of $H^\sigma(K)$. The proof is essential the same as that of Theorem \ref{thm69}, which relies on Lemma \ref{lemma68}. Below we provide a lemma analogous to Lemma \ref{lemma68}, and omit the proof of Theorem \ref{thm46} (a).

\begin{lemma}\label{lemma75}
	(a). For any $m'\geq m\geq 1$ and any $f\in \mathcal{H}_0\oplus(\oplus_{l=1}^m T_l)$, we have
	\[\|\nabla_{m'} f\|_{l^2(E_{v,m'})}\lesssim r^{\frac{m'-m}{2}}\|\nabla_m f\|_{l^2(E_{v,m})}.\]
	
	(b). For any $m'> m\geq 1$ and any $f\in T_{m'}$, we have $\nabla_{m}f=0$.
\end{lemma}
\textit{Proof.} (a). We have the estimate that
\[\mathcal{E}^{1/2}(f)\asymp r^{-m'/2}\|\nabla_{m'}f\|_{l^2(E_{v,m'})},\]
since $f$ is harmonic in $F_wK$ for any $w\in \Lambda_m$. Thus,
\[\|\nabla_{m'}f\|_{l^2(E_{v,m'})}\asymp r^{m'/2}\mathcal{E}^{1/2}(f)\asymp r^{(m'-m)/2}\|\nabla_m f\|_{l^2(E_{v,m})}.\]
(b) is trivially true since $f|_{F_wV_0}=0$ for any $w\in \Lambda_{m'-1}$.\hfill$\square$

\subsection{Proof of  $\Gamma^\sigma(K)=\tilde{\Gamma}^\sigma(K)\cap C(K)$ for $\sigma>\frac{d_S}{2}$}
In this part, we prove Theorem \ref{thm59}, which illustrates the relation between $\Gamma^\sigma(K)$ and $\tilde{\Gamma}^\sigma(K)$ when $\sigma>\frac{d_S}{2}$. Since $H^\sigma(K)=\Gamma^\sigma(K)$ when $\sigma<1$, this gives another characterization of $H^\sigma(K)$.

First, we introduce the following notations.

For each $\omega\in \pi^{-1}(V_*)$, define
\[\{w(k,\omega)\}_{k=0}^\infty=\big\{[\omega]_l:l\geq 0\text{, }[\omega]_l\in\bigcup_{m=0}^\infty \Lambda_m\text{ and }\pi(\omega)\in F_{[\omega]_l}V_0\big\},\]
with order \[|w(0,\omega)|<|w(1,\omega)|<|w(2,\omega)|<\cdots.\] We then easily have the following properties.

\textit{1. If $\pi(\omega)\in \big(\bigcup_{w\in \Lambda_m} F_wV_0\big)\setminus \big(\bigcup_{w\in \Lambda_{m-1}} F_wV_0\big)$, then $w(k,\omega)\in\Lambda_{m+k}$;}

\textit{2. Let $w\in\Lambda_m$ and $\omega\in \pi^{-1}(V_*)$, we have $w\in \{w(k,\omega)\}_{k=0}^\infty$ if and only if $[\omega]_{|w|}=w$ and $\sigma^{|w|}(\omega)\in\mathcal{P}$};

\textit{3. $\bigcup_{\omega\in \pi^{-1}(V_*)}\{w{(k,\omega)}\}_{k=0}^\infty=\bigcup_{m=0}^\infty \Lambda_m.$}

We have the following estimates.
\begin{lemma}\label{lemma76} Let $\sigma\geq 0$. For $f=C+\sum_{m=1}^\infty\tilde{f}_m\in \tilde{\Gamma}^\sigma(K)$ with $C\in\mathbb{R}$ and $\tilde{f}_m\in \tilde{J}_m$, we have
	\[\sum_{\omega\in \pi^{-1}(V_*)}\sum_{k=0}^\infty r_{w{(k,\omega)}}^{d_H-\sigma d_W}\big(A_{w{(k+1,\omega)}}(f)-A_{w{(k,\omega)}}(f)\big)^2\lesssim\sum_{m=1}^\infty r^{-md_H}\lambda^{2m}\|\tilde{f}_m\|^2_{L^2(K)}.\]
\end{lemma}
\textit{Proof.} Recall $W_w=\{w'\in W_*:F_{w'}K\subset F_wK\}$ in Definition \ref{def56}. It is direct to check that
\[\|\tilde{f}_m\|^2_{L^2(K)}=\sum_{w\in \Lambda_{m-1}}\sum_{w'\in W_w\cap \Lambda_m}r^{d_H}_{w'}\big(A_{w'}(f)-A_w(f)\big)^2,\]
and thus
\[\begin{aligned}
\sum_{m=1}^\infty r^{-md_H}\lambda^{2m}\|\tilde{f}_m\|^2_{L^2(K)}&\asymp \sum_{m=1}^\infty\lambda^{2m}\sum_{w\in \Lambda_{m-1}}\sum_{w'\in W_w\cap \Lambda_m}\big(A_{w'}(f)-A_w(f)\big)^2\\
&\gtrsim \sum_{\omega\in \pi^{-1}(V_*)}\sum_{k=0}^\infty r_{w{(k,\omega)}}^{d_H-\sigma d_W}\big(A_{w{(k+1,\omega)}}(f)-A_{w{(k,\omega)}}(f)\big)^2,
\end{aligned}\]
since each $\big(A_{w'}(f)-A_w(f)\big)^2$ occurs at most $\#\mathcal{P}$ times in the last term of the above estimate by Property 2.\hfill$\square$\vspace{0.2cm}

\begin{lemma}\label{lemma77}
	Let $\sigma>\frac{d_S}{2}$. For $f\in\tilde{\Gamma}^\sigma(K)\cap C(K)$, we have
	\[\sum_{\omega\in \pi^{-1}(V_*)}\sum_{k=0}^\infty r_{w{(k,\omega)}}^{d_H-\sigma d_W}\big(A_{w{(k,\omega)}}(f)-f(\pi(\omega))\big)^2\lesssim \|f\|^2_{\tilde{\Gamma}^\sigma(K)}.\]
\end{lemma}

\textit{Proof.} Fix $\omega\in \pi^{-1}(V_*)$. By Property 1, for any $k\geq 0$, we have $w{(k,\omega)}\in \Lambda_{k+m}$ for some $m\geq 0$. Noticing that $f(\pi(\omega))=\lim_{k\to\infty}A_{w(k,\omega)}(f)$, we have
\[\begin{aligned}
\big\|r_{w{(k,\omega)}}^{(d_H-\sigma d_W)/2}\big(A_{w{(k,\omega)}}(f)-f(\pi(\omega))\big)\big\|_{l^2}&\asymp \lambda^m\big\|\lambda^k\big(A_{w{(k,\omega)}}(f)-f(\pi(\omega))\big)\big\|_{l^2}\\
&=\lambda^m\big\|\sum_{l=0}^\infty\lambda^{-l}\lambda^{k+l}\big(A_{w_{(k+l,\omega)}}(f)-A_{w_{(k+l+1,\omega)}}(f)\big)\big\|_{l^2}\\
&\lesssim \lambda^m\big\|\lambda^k\big(A_{w(k,\omega)}(f)-A_{w(k+1,\omega)}(f)\big)\big\|_{l^2}\\
&\asymp \big\|r_{w{(k,\omega)}}^{(d_H-\sigma d_W)/2}\big(A_{w(k+1,\omega)}(f)-A_{w(k,\omega)}(f)\big)\big\|_{l^2},
\end{aligned}\]
by using Minkowski inequality and the fact that $\lambda>1$.  Thus,
\[\begin{aligned}
&\sum_{\omega\in \pi^{-1}(V_*)}\sum_{k=0}^\infty r_{w{(k,\omega)}}^{d_H-\sigma d_W}\big(A_{w{(k,\omega)}}(f)-f(\pi(\omega))\big)^2\\\lesssim &\sum_{\omega\in \pi^{-1}(V_*)}\sum_{k=0}^\infty r_{w{(k,\omega)}}^{d_H-\sigma d_W}\big(A_{w{(k+1,\omega)}}(f)-A_{w{(k,\omega)}}(f)\big)^2.
\end{aligned}\]
The lemma follows immediately from Lemma \ref{lemma76}. \hfill$\square$\vspace{0.2cm}

\begin{lemma}\label{lemma78}
	Let $\sigma>\frac{d_S}{2}$. For $f\in\tilde{\Gamma}^\sigma(K)\cap C(K)$, we have $f\in\Lambda^{2,2}_{\sigma}(K)$ with $\|f\|_{\Lambda^{2,2}_{\sigma}(K)}\lesssim \|f\|_{\tilde{\Gamma}^\sigma(K)}$.
\end{lemma}
\textit{Proof.} The lemma is true since
\[\begin{aligned}
\sum_{m=0}^\infty \lambda^{2m}\|\nabla_m f\|^2_{l^2(E_{v,m})}&\asymp\sum_{m=0}^\infty\lambda^{2m}\sum_{w\in \Lambda_m}\sum_{p\neq q\in V_0}\big(f(F_wp)-f(F_wq)\big)^2\\
&\lesssim\sum_{m=0}^\infty \sum_{w\in \Lambda_m}r_w^{d_H-\sigma d_W}\sum_{p\in V_0}\big(A_w(f)-f(F_wp)\big)^2\\
&\lesssim\sum_{\omega\in \pi^{-1}(V_*)}\sum_{k=0}^\infty r_{w{(k,\omega)}}^{d_H-\sigma d_W}\big(A_{w{(k,\omega)}}(f)-f(\pi(\omega))\big)^2\\
&\lesssim \|f\|^2_{\tilde{\Gamma}^\sigma(K)},
\end{aligned}\]
where we use Lemma \ref{lemma77} in the last line.\hfill$\square$\vspace{0.2cm}

\textit{Proof of Theorem \ref{thm59}}.  It is obvious that $\Gamma^\sigma(K)\subset\tilde{\Gamma}^\sigma(K)$ by Proposition \ref{prop57}. In addition, we can easily see $\Gamma^\sigma(K)\subset C(K),\forall \sigma>\frac{d_S}{2}$ by Theorem \ref{thm69} and $\|S_l\tilde{f}_l\|_{L^\infty(K)}\lesssim\|D_lS_l\tilde{f}_l\|_{l^2(E_{c,l})},\forall \tilde{f}_l\in\tilde{J}_l$. Thus, we have $\Gamma^\sigma(K)\subset\tilde{\Gamma}^\sigma(K)\cap C(K)$.

Conversely, by Lemma \ref{lemma78}, we have $\tilde{\Gamma}^\sigma(K)\cap C(K)\subset\Lambda^{2,2}_\sigma(K)$ whenever $\sigma>\frac{d_S}{2}$. For $\frac{d_S}{2}<\sigma<1$, this gives that $\tilde{\Gamma}^\sigma(K)\cap C(K)\subset\Gamma^\sigma(K)$ by Theorem \ref{thm46} (a) and \ref{thm61}. For $\sigma\geq 1$, by the Remark after Definition \ref{def43}, we then have $\tilde{\Gamma}^\sigma(K)\cap C(K)\subset constants$. The desired result follows.
\hfill$\square$

\section{Atomic decompositions: higher orders}
In this section, we discuss the proof of the atomic decomposition theorem (Theorem \ref{decompositionthm}) of the spaces $H^\sigma(K)$ for $1\leq\sigma<2$. Since the story is parallel to the case $0\leq \sigma<1$, we will omit some details in the proof.  As a byproduct, we will prove Theorem \ref{thm46} (b).

As in previous sections, we write $\lambda=\lambda(\sigma)=r^{(d_H-\sigma d_W)/2}$.

\subsection{Proof of  $H^\sigma(K)=\tilde{\Lambda}^{2,2}_{\sigma}(K)$ with $\frac{d_S}{2}<\sigma<2$} In this part, we develop necessary lemmas for the atomic decomposition for $1\leq\sigma<2$, analogous to Section 5 and 6. As an application, we will prove Theorem \ref{thm46} (b).

Recall $H_{\Lambda_m}$ defined in Definition \ref{def44}. We have some easy estimates for tent functions concerning $H_{\Lambda_m}$.

\begin{lemma}\label{lemma81}
	For $m\geq l$ and $\varphi_l\in T_l$, we have $$\|H_{\Lambda_m}\varphi_l\|_{l^2(V_{\Lambda_m}\setminus V_0)}\asymp r^{-l}\|\nabla_l \varphi_l\|_{l^2(E_{v,l})}.$$
\end{lemma}
\textit{Proof.} First, we take $m=l$. It is not hard to see that
\[\|H_{\Lambda_l}\varphi_l\|_{l^2(V_{\Lambda_l}\setminus V_0)}\asymp r^{-l} \|\nabla_l \varphi_l\|_{l^2(E_{v,l})}.\]
Next, we claim that $\|H_{\Lambda_m}\varphi_l\|_{l^2(V_{\Lambda_m}\setminus V_0)}=\|H_{\Lambda_l}\varphi_l\|_{l^2(V_{\Lambda_l}\setminus V_0)}$ for $m\geq l$. To see this, we notice that the following equality holds for any $g\in dom\mathcal{E}$,
\[-(g,H_{\Lambda_l}\varphi_l)=\mathcal{E}_{\Lambda_l}(g,\varphi_l)=\mathcal{E}(g,\varphi_l)=\mathcal{E}_{\Lambda_m}(g,\varphi_l)=-(g,H_{\Lambda_m}\varphi_l).\]
This implies $H_{\Lambda_l}\varphi_l(x)=H_{\Lambda_m}\varphi_l(x),\forall x\in V_{\Lambda_l}$ and $H_{\Lambda_m}\varphi_l(x)=0,\forall x\in V_{\Lambda_m}\setminus V_{\Lambda_l}$, as $g$ can be arbitrarily chosen. \hfill$\square$\vspace{0.2cm}

Using Lemma \ref{lemma81}, by a similar proof of Theorem \ref{thm58}, we have

\begin{lemma}\label{lemma82}
	Let $\frac{d_S}{2}<\sigma<
	2-\frac{d_S}{2}$ and $f\in C(K)$. Write $f=\sum_{m=0}^\infty\varphi_m$ with $\varphi_0\in\mathcal{H}_0$ and $\varphi_m\in T_m,\forall m\geq 1$. Then we have
	\[\begin{aligned}
	\|\varphi_0\|^2_{L^2(K)}+\sum_{m=0}^\infty \lambda^{2m}\|\nabla_m\varphi_m\|^2_{l^2(E_{v,m})}
	\asymp \|\varphi_0\|^2_{L^2(K)}+\sum_{m=1}^\infty {r^{2m}}\lambda^{2m}\|H_{\Lambda_m}f\|^2_{l^2(V_{\Lambda_m}\setminus V_0)}.
	\end{aligned}\]
\end{lemma}

Next, we introduce the following \textit{smoothed tent functions}, analogous to smoothed Haar functions.

\begin{definition}\label{def83}
	For $m\geq 1$, for any $\varphi_m\in T_m$, define $\breve{S}_m\varphi_m$ to be the unique function in $H^2(K)$ such that
	\[\breve{S}_m\varphi_m(x)=\varphi_m(x), \forall x\in \bigcup_{w\in \Lambda_m}F_wV_0,\]
	and
	\[\|\Delta \breve{S}_m\varphi_m\|_{L^2(K)}=\min\{\|\Delta f\|_{L^2(K)}:f\in H^2(K), f|_{V_{\Lambda_m}}=\varphi_m|_{V_{\Lambda_m}}\},\]
	call it a $m$-smoothed tent function.
	
	Define the space of $m$-smoothed tent functions by $\breve{T}_m=\breve{S}_mT_m$.
\end{definition}

In particular, we have the decomposition theorem (Theorem \ref{thm84}) analogous to Theorem \ref{prop67} and \ref{thm69}.

\begin{theorem}\label{thm84}
	For $\frac{d_S}{2}<\sigma\leq 2$, each function $f$ in $H^\sigma(K)$ admits a unique series expansion $f=\sum_{m=0}^\infty \breve{\varphi}_m$ with $\breve{\varphi}_0\in\mathcal{H}_0$ and $\breve{\varphi}_m\in \breve{T}_m$, and
	\begin{equation}
	\|f\|_{H^\sigma(K)}\asymp \big(\|\breve{\varphi}_0\|^2_{L^2(K)}+\sum_{m=0}^\infty\lambda^{2m}\|\nabla_m\breve{\varphi}_m\|^2_{l^2(E_{v,m})}\big)^{1/2}.
	\end{equation}
	Moreover, for any series $f=\sum_{m=0}^\infty \breve{\varphi}_m$, $f\in H^\sigma(K)$ if and only if $$\sum_{m=0}^\infty\lambda^{2m}\|\nabla_m\breve{\varphi}_m\|^2_{l^2(E_{v,m})}<\infty.$$
\end{theorem}

The proof of Theorem \ref{thm84} is similar to that in Section 6. Below we list a key lemma but omit its proof.

\begin{lemma}\label{lemma85}
	(a). For $m\neq m'\geq 1$, we have $\Delta (\breve{T}_m)\bot \Delta (\breve{T}_{m'})$ in $L^2(K)$.
	
	(b). For  $m\geq 1$ and $f\in\mathcal{H}_0\oplus(\oplus_{l=1}^m\breve{T}_l)$, we have
	\[\|\Delta f\|_{L^2(K)}\asymp r^{-md_H/2}\|H_{\Lambda_m}f \|_{l^2(V_{\Lambda_m}\setminus V_0)}.\]
	In particular, for any $m\geq 1$ and $\breve{\varphi}_m\in\breve{T}_m$, we have $$\|\Delta\breve{\varphi}_m\|_{L^2(K)}\asymp r^{-md_H/2}\|H_{\Lambda_m}\breve{\varphi}_m \|_{l^2(V_{\Lambda_m}\setminus V_0)}\asymp r^{-m(1+d_H/2)}\|\nabla_m\breve{\varphi}_m\|_{l^2(E_{v,m})}.$$
	
	(c). For $m\geq 1$ and  $f\in H^2(K)$, we have
	\[\|H_{\Lambda_m} f\|_{l^2(V_{\Lambda_m}\setminus V_0)}\lesssim r^{md_H/2}\|\Delta f\|_{L^2(K)}.\]
	In particular, for $m\geq l$ and  $\breve{\varphi}_l\in\breve{T}_l$, we have
	\[\|H_{\Lambda_m} \breve{\varphi}_l\|_{l^2(V_{\Lambda_m}\setminus V_0)}\lesssim r^{md_H/2}\|\Delta \breve{\varphi}_l\|_{L^2(K)}\asymp r^{md_H/2-l(1+d_H/2)}\|\nabla_l \breve{\varphi}_l\|_{l^2({E_{v,l}})}.\]
\end{lemma}

Using Lemma \ref{lemma85} (a) and (b), we get the decomposition of $H^2(K)$ analogous to Theorem \ref{prop67}. Using Lemma \ref{lemma85} (b) and (c), combining with Lemma \ref{lemma82} and the atomic decomposition Theorem for $\frac{d_S}{2}<\sigma<1$, we get the decomposition of $H^\sigma(K)$ for $\frac{d_S}{2}<\sigma<1$ analogous to Theorem \ref{thm69}. Then we finish the proof of Theorem \ref{thm84} by using complex interpolation.

\vspace{0.2cm}
\textit{Proof of Theorem \ref{thm46} (b).}  It follows easily from Theorem \ref{thm84} and Lemma \ref{lemma85} (c). \hfill$\square$
\vspace{0.2cm}

\subsection{Atomic decomposition: $1\leq\sigma<2$}
The case $1\leq\sigma<2-\frac{d_S}{2}$  follows immediately from Lemma \ref{lemma72}, Lemma \ref{lemma82} and Theorem \ref{thm46} (b). It remains to consider the case $2-\frac{d_S}{2}<\sigma<2$.

The proof relies on the following lemmas. Note that by Theorem \ref{thm84}, for $2-\frac{d_S}{2}<\sigma<2$, each $f\in H^\sigma(K)$ admits a unique expansion $f=\sum_{m=0}^{\infty}\breve{\varphi}_m$ with $\breve{\varphi}_0\in\mathcal{H}_0$ and $\breve{\varphi}_m\in\breve{T}_m$, $\forall m\geq 1$.

\begin{lemma}\label{lemma86} Let $l\geq 1$, $\breve{\varphi}_l=GC_l+\sum_{m=1}^\infty G\tilde{f}_{m,l}$ with $C_l\in \mathbb{R}$ and $\tilde{f}_{m,l}\in \tilde{J}_m$.
	
	(a). If $m<  l$, then $\|\tilde{f}_{m,l}\|_{L^2(K)}\lesssim r^{(l-m)d_H/2-l(1+d_H/2)}\|\nabla_l\breve{\varphi}_l\|_{l^2(E_{v,l})}$.
	
	(b). If $m\geq l$, then $\|\tilde{f}_{m,l}\|_{L^2(K)}\lesssim r^{-l(1+d_H/2)}\|\nabla_l\breve{\varphi}_l\|_{l^2(E_{v,l})}$.
	
	In particular, $|C_l|\lesssim r^{-l}\|\nabla_l\breve{\varphi}_l\|_{l^2(E_{v,l})}$.
\end{lemma}
\textit{Proof.} Note that $-\Delta\breve{\varphi}_l=C_l+\sum_{m=1}^\infty\tilde{f}_{m,l}\in L^2(K)$, and by Lemma \ref{lemma85}, $\|\Delta\breve{\varphi}_l\|_{L^2(K)}\asymp r^{-l(1+d_H/2)}\|\nabla_l\breve{\varphi}_l\|_{l^2(E_{v,l})}$. Obviously, (b) is trivial. To prove (a),  fix $w\in \Lambda_m$, and define a function $\psi_w$ which takes value
\[\psi_w(x)=\begin{cases}
1,\text{ if }x\in  V_{\Lambda_{l-1}}\cap(F_wK\setminus F_wV_0),\\
0,\text{ if }x\in V_{\Lambda_{l-1}}\setminus (F_wK\setminus F_wV_0),
\end{cases}\]
and is harmonic in each $F_{w'}K,w'\in \Lambda_{l-1}$. Then we have
\[\int_{F_wK}\Delta\breve{\varphi}_l\cdot\psi_w d\mu=\int_{K}\Delta\breve{\varphi}_l\cdot\psi_w d\mu =-\mathcal{E}(\breve{\varphi}_l,\psi_w)=-\mathcal{E}_{\Lambda_{l-1}}(\breve{\varphi}_l,\psi_w)=0.\]
As a consequence, we have the estimate that
\[\begin{aligned}
|A_w(\Delta\breve{\varphi}_l)|&=\big|r_w^{-d_H}\int_{F_wK}\Delta\breve{\varphi}_ld\mu\big|=\big|r_w^{-d_H}\int_{F_wK}(1-\psi_w)\Delta\breve{\varphi}_ld\mu\big|\\&\lesssim r^{-md_H}r^{ld_H/2}\|\Delta\breve{\varphi}_l\|_{L^2(F_wK)},
\end{aligned}\]
and thus
\[r^{md_H}|A_w(\Delta\breve{\varphi}_l)|^2\lesssim r^{(l-m)d_H}\|\Delta\breve{\varphi}_l\|^2_{L^2(F_wK)},
\]
Thus (a) follows by summing up the above estimates over all cells $F_wK$ with $w\in \Lambda_m$ and using the fact $\|\Delta\breve{\varphi}_l\|_{L^2(K)}\asymp r^{-l(1+d_H/2)}\|\nabla_l\breve{\varphi}_l\|_{l^2(E_{v,l})}$. The estimate of $|C_l|$ follows by a same argument as that of (b).\hfill$\square$

\begin{lemma}\label{lemma87}
	Let $l\geq 1$, $\tilde{f}_l\in \tilde{J}_l$, write $G\tilde{f}_l=\sum_{m=1}^\infty \breve{\varphi}_{m,l}$ with  $\breve{\varphi}_{m,l}\in \breve{T}_m$.
	
	(a). If $m<l$, then $\|\nabla_m\breve{\varphi}_{m,l}\|_{l^2(E_{v,m})}\lesssim r^{(l-m) d_W/2+m(1+d_H/2)}\|\tilde{f}_l\|_{L^2(K)}$.
	
	(b). If $m\geq l$, then  $\|\nabla_m\breve{\varphi}_{m,l}\|_{l^2(E_{v,m})}\lesssim r^{m(1+d_H/2)}\|\tilde{f}_l\|_{L^2(K)}$.
\end{lemma}

\textit{Proof.} By Theorem \ref{thm84}, the expansion $G\tilde{f}_l=\sum_{m=1}^\infty \breve{\varphi}_{m,l}$ follows with
\[\|G\tilde{f}_l\|_{H^2(K)}\asymp \big(\sum_{m=1}^\infty r^{-2m(1+d_H/2)}\|\nabla_m\breve{\varphi}_{m,l}\|_{l^2(E_{v,m})}^2\big)^{1/2}.\]
This trivially gives (b). To prove (a), we first estimate $H_{\Lambda_m}G\tilde{f}_l$ for a fixed $m<l$. For each $x\in V_{\Lambda_m}\setminus V_0$, we have
\[H_{\Lambda_m}G\tilde{f}_l(x)=-\mathcal{E}(\psi_x^{(m)},G\tilde{f}_l)=\int_K\psi_x^{(m)}\tilde{f}_ld\mu,\]
where $\psi^{(m)}_x$ is a piecewise harmonic function which takes value
\[
\psi^{(m)}_x(y)=\begin{cases}1,\text{ if }y=x,\\
0,\text{ if }x\in V_{\Lambda_m}\setminus \{x\},\end{cases}
\]
and is harmonic in each $F_wK$ with $w\in\Lambda_m$. Note that $\psi^{(m)}_x=\psi_x$ if $x\in V_{\Lambda_m}\setminus V_{\Lambda_{m-1}}$ in our previous notation for atomic decomposition theory. Denote  $P_{\tilde{J}_l}$ the orthogonal projection from $L^2(K)$ onto $\tilde{J}_l$. We see that
\[\|P_{\tilde{J}_l}\psi^{(m)}_x\|_{L^2(K)}\asymp r^{ld_H/2}\|D_lP_{\tilde{J}_l}\psi^{(m)}_x\|_{l^2(E_{c,l})}\lesssim r^{ld_H/2}r^{l/2}\mathcal{E}^{1/2}(\psi^{(m)}_x)\lesssim r^{ld_H/2+(l-m)/2}.\]
Thus
\[\begin{aligned}
\big|H_{\Lambda_m}G\tilde{f}_l(x)\big|&=\big|\int_K\psi^{(m)}_x\tilde{f}_ld\mu\big|=\big|\int_K P_{\tilde{J}_l}\psi^{(m)}_x\cdot \tilde{f}_ld\mu\big|\\&\lesssim r^{ld_H/2+(l-m)/2}\|\tilde{f}_l\|_{L^2(supp P_{\tilde{J}_l}\psi^{(m)}_x)}.
\end{aligned}\]
Note that $P_{\tilde{J}_l}\psi^{(m)}_x$ is locally supported in $\bigcup\{F_wK:x\in F_wK, w\in\Lambda_m\}$. Summing the above estimate over vertices in $V_{\Lambda_m}\setminus V_0$, we get
\[\|H_{\Lambda_m}\sum_{m'=1}^m \breve{\varphi}_{m',l}\|_{l^2(V_{\Lambda_m}\setminus V_0)}=\|H_{\Lambda_m}G\tilde{f}_l\|_{l^2(V_{\Lambda_m}\setminus V_0)}\lesssim r^{ld_H/2+(l-m)/2}\|\tilde{f}_l\|_{L^2(K)}.\]
Then by Lemma \ref{lemma85} (a) and (b), we get
\[\big\|\Delta\breve{\varphi}_{m,l}\big\|_{L^2(K)}\lesssim \|\Delta\sum_{m'=1}^m \breve{\varphi}_{m',l}\|_{L^2(K)}\lesssim r^{(l-m) d_W/2}\|\tilde{f}_l\|_{L^2(K)},\]
and thus by Lemma \ref{lemma85} (c),
(a) follows immediately. \hfill$\square$\\

The proof of Theorem \ref{decompositionthm} for $2-\frac{d_S}{2}<\sigma<2$ is essentially the same as case for $\frac{d_S}{2}<\sigma<1$ in Section 7, by using Lemma \ref{lemma86} and \ref{lemma87} instead of Lemma \ref{lemma73} and \ref{lemma74}. Readers only need to carefully check the orders involved. We omit it.

\bibliographystyle{amsplain}

\end{document}